\documentclass[12pt,a4paper,leqno]{amsart} %Fr o m detta
\usepackage{amssymb}
\usepackage{amsthm}
\usepackage{latexsym}
\usepackage{amsmath}
\usepackage{mathrsfs}
\usepackage[X2,T1]{fontenc}
\usepackage{calc}                         %Detta behövs för att fixa
\usepackage{cite}

\newcommand{\scal}[2]{\langle #1,#2\rangle}

\newcommand{\rr}[1]{\mathbf R^{#1}}
\newcommand{\co}{\mathbf C}

\newcommand{\nm}[2]{\Vert #1\Vert _{#2}}
\newcommand{\nmm}[1]{\Vert #1\Vert }
\newcommand{\abp}[1]{\vert #1\vert}

\newcommand{\op}{\operatorname{Op}}

\newcommand{\sets}[2]{\{ \, #1\, ;\, #2\, \} }
\newcommand{\ep}{\varepsilon}
\newcommand{\fy}{\varphi}

\newcommand{\cdo}{\, \cdot \, }

\newcommand{\wpr}{{\text{\footnotesize $\#$}}}

\newcommand{\eabs}[1]{\langle #1\rangle}
\newcommand{\ON}{\operatorname{ON}}

\newcommand{\tp}{\operatorname{Tp}}
\newcommand{\vrum}{\vspace{0.1cm}}

\numberwithin{equation}{section}          %Detta gör att man får
\newtheorem{thm}{Theorem}
\numberwithin{thm}{section}
\newtheorem*{tom}{\rubrik}
\newcommand{\rubrik}{}
\newtheorem{prop}[thm]{Proposition}
\newtheorem{cor}[thm]{Corollary}
\newtheorem{lemma}[thm]{Lemma}

\theoremstyle{definition}

\newtheorem{defn}[thm]{Definition}

\theoremstyle{remark}

\newtheorem{rem}[thm]{Remark}              %T o m hit är bara allmän
\author{Joachim Toft}

\address{Department of Mathematics and Systems Engineering,
V{\"a}xj{\"o} University, Sweden}

\email{joachim.toft@vxu.se}

\title{Products on Schatten-von Neumann classes and modulation spaces}


\subjclass[2000]{??}

\begin{document}

\begin{abstract}
We consider modulation space and spaces of Schatten-von Neumann
symbols where corresponding pseudo-differential operators map
one Hilbert space to another. We prove 
H{\"o}lder-Young and Young type results for such spaces under dilated
convolutions and multiplications. We also prove continuity properties
for such spaces under the twisted convolution, and the Weyl
product. These results lead to continuity
properties for twisted convolutions on Lebesgue spaces,
e.{\,}g. $L^p_{(\omega )}$ is a twisted convolution algebra when
$1\le p\le 2$ and appropriate weight $\omega$.
\end{abstract}

\maketitle

\vspace{0.5cm}

\par

\section{Introduction}\label{sec0}

\par

In this paper we establish continuity properties for various
products on modulation spaces and a familly of symbol classes such that
corresponding pseudo-differential operators are of
Schatten-von Neumann types.

 This means that each symbol class consists of all
tempered distributions such that the corresponding pseudo-differential
operators are Schatten-von Neumann operators of certain degree from
one Hilbert space to another.
For such spaces of functions and distributions, we
establish Young type and H{\"o}lder-Young type inequalities with
respect to the Weyl product, twisted convolution, dilated convolutions
and dilated multiplications. These products are important in the
theory of pseudo-differential operators. In fact, the Weyl product
corresponds to compositions of Weyl operators on the symbol side. On
the symplectic Fourier transform side, the Weyl product takes the form
as a twisted convolution. In the theory of pseudo-differential
operators, it is in many situations convenient to approximate a
pseudo-differential operator with a Toeplitz operator. Then the Weyl
symbol of a Toeplitz operator is an ordinary convolution of the
Toeplitz symbol and a rank one element, which can be rewritten as a
product between symbols by using the symplectic Fourier transform.

\par

In particular we generalize results in
\cite{To3,To5}, where similar questions were considered for classical
modulation spaces, and spaces of pseudo-differential operators of
Schatten-von Neumann types on $L^2$. The Weyl product in the context
of modulation space theory was recently investigated in
\cite{HTW}. From these results we establish continuity for the twisted
convolutions when acting on modulation spaces of Wiener amalgam type, using
the fact that the Fourier transform of a Weyl product
is essentially a twisted convolution of the Weyl symbols. From these
results we thereafter prove continuity properties of
the twisted convolution on weighted Fourier Lebesgue spaces.

\par

We also apply our results on Toeplitz operators and prove that for
each appropriate dilated symbol to a Schatten-von Neumann
pseudo-differential operator, then the Toeplitz operator belongs the
same Schatten-von Neumann class.

\par

In order to be more specific, we recall some definitions. Assume that
$t\in \mathbf R$ is fixed and that $a\in \mathscr S(\rr {2d})$. (We
use the same notation for the  usual function and distribution spaces
as in \cite{Ho1}.) Then the pseudo-differential operator $a_t(x,D)$
with symbol $a$ is the continuous operator on $\mathscr S(\rr d)$,
defined by the formula
\begin{equation}\label{e0.5}
\begin{aligned}
(a_t(x,D)f)(x &)
=
(\op _t(a)f)(x)
\\
=
(2\pi  &) ^{-d}\iint a((1-t)x+ty,\xi )f(y)e^{i\scal {x-y}\xi }\,
dyd\xi .
\end{aligned}
\end{equation}
The definition of $a_t(x,D)$ extends to each $a\in \mathscr S'(\rr
{2d})$, and then $a_t(x,D)$ is continuous from $\mathscr S(\rr d)$ to
$\mathscr S'(\rr d)$. (Cf. e.{\,}g. \cite {Ho1}.) In fact, the Fourier
transform $\mathscr F$ is the linear and continuous operator on
$\mathscr S'(\rr d)$ which takes the form
\begin{equation}\label{ftransf}
(\mathscr Ff)(\xi ) =\widehat f(\xi )=(2\pi )^{-d/2}\int f(x)e^{i\scal
x\xi}\, dx,
\end{equation}
when $f\in L^1(\rr d)$. If $\mathscr F_2F$ is the partial Fourier
transform of $F(x,y)\in \mathscr S'(\rr {2d})$ with respect to the $y$
variable and $a\in \mathscr S'(\rr {2d})$, then we let
$a_t(x,D)$ be the linear and continuous operator from $\mathscr S(\rr
d)$ to $\mathscr S'(\rr d)$ with the distribution kernel
\begin{equation}\label{atkernel}
K_{a,t}(x,y)=(\mathscr F_2^{-1}a)((1-t)x+ty,x-y).
\end{equation}

\par

If $t=1/2$, then $a_t(x,D)$ is equal to the Weyl operator $a^w(x,D)$
for $a$. If instead $t=0$, then the standard (Kohn-Nirenberg)
representation $a(x,D)$ is obtained.

\par

For each $a\in \mathscr S(\rr {2d})$ we also let $Aa$ be the linear
and continuous operator from $\mathscr S(\rr d)$ to $\mathscr
S(\rr d)$, given by
\begin{equation}\label{Aadef}
(Aa)(x,y) = (2\pi )^{-d/2}\int a((y-x)/2,\xi )e^{-i\scal {x+y}\xi}\,
d\xi .
\end{equation}
(Here and in what follows we identify operators with their
distribution kernels.)
This operator representation is closely related to the Weyl
quantization in the sense that
$$
Aa(x,y)=(2\pi )^{d/2}K_{a,1/2}(-x,y) = (2\pi )^{d/2}K_{\mathscr
F_\sigma a,1/2}(x,y),
$$
where $\mathscr F_\sigma$ denotes the symplectic Fourier transform
(see Section \ref{sec1} for strict definitions). Here the first
equality follows from \eqref{atkernel} and \eqref{Aadef}. In
particular, the definition of $Aa$ extends to each $a\in \mathscr
S'(\rr {2d})$ and then $Aa$ is continuous from $\mathscr S(\rr d)$ to
$\mathscr S'(\rr d)$, since the same is true for the Weyl
quantization.

\par

In the Weyl calculus, operator composition corresponds on the symbol
level to the \emph{Weyl product}, sometimes also called the twisted
product, denoted by $\wpr$. In Sections \ref{sec2}--\ref{sec3} we
establish extensions for the Weyl product and twisted convolutions on
modulation spaces, which we shall describe now.

\par

The modulation spaces was introduced by Feichtinger in \cite{Fe4}, and
developed further and generalized in \cite{Fe5,FG1,FG2,FG3,Gc2}, where
Feichtinger and Gr{\"o}chenig established the theory of coorbit
spaces. In particular, the modulation spaces $M^{p,q}_{(\omega )}$ and
$W^{p,q}_{(\omega )}$, where $\omega$ denotes a weight function on
phase (or time-frequency shift) space, appear as the set of tempered
(ultra-) distributions
whose STFT belong to the weighted and mixed Lebesgue space
$L^{p,q}_{1,(\omega )}$ and $L^{p,q}_{2,(\omega )}$ respectively. (See
Section \ref{sec1} for strict definitions.) It follows that $\omega$,
$p$ and $q$ to some extent quantify the degrees of
asymptotic decay and singularity of the distributions in
$M^{p,q}_{(\omega)}$ and $W^{p,q}_{(\omega )}$. By choosing the weight
$\omega$ in appropriate ways, the space $W^{p,q}_{(\omega )}$ becomes
a Wiener amalgam space, introduced by Feichtinger in
\cite{Fe2}. Furthermore, if $\omega$ is trivial, i.{\,}e. $\omega =1$,
then $M^{p,q}_{(\omega)}$ is the \emph{classical} modulation space
$M^{p,q}$. (See \cite{Fe6} for the most updated description of
modulation spaces.)

\par

From the construction of these spaces, it turns out that modulation
spaces of the form $M^{p,q}_{(\omega )}$ and Besov spaces in some
sense are rather similar, and sharp embeddings between these
spaces can be found in \cite{To5, To7}, which
are improvements of certain embeddings in \cite {Gb}. (See
also \cite {ST1} for verification of the sharpness.) In the same
way it follows that modulation spaces of the form $W^{p,q}_{(\omega
)}$ and Triebel-Lizorkin spaces are rather similar.

\par

During the last 15 years many results have been proved which confirm
the usefulness of the modulation
spaces and their Fourier transforms in time-frequency analysis, where
they occur naturally. For example, in
\cite{FG3,Gc2}, it is shown that all
such spaces admit reconstructible sequence space representations
using Gabor frames.

\par

Parallel to this development, modulation spaces have been incorporated
into the calculus of pseudo-differential operators. In fact, in
\cite{GH1,Gc2}, Gr{\"o}chenig and Heil prove that each operator with
symbol in $M^{\infty,1}$ is continuous on all modulation spaces
$M^{p,q}$, $p,q \in [1,\infty]$. In particular, since $M^{2,2}=L^2$ it
follows that any such operator is $L^2$ continuous, which was proved
by Sj{\"o}strand already in \cite{Sj1}.

Some further generalizations to operators with
symbols in more general modulation spaces were obtained in
\cite{GH2,To5, To6, To8}. Modulation spaces in
pseudodifferential calculus is currently an active field of research
(see e.{\,}g.
\cite{GH2, HRT, La1, ST1, Ta, Te1, To2, To5, To8}).

\par

In the Weyl calculus of pseudo-differential operators, operator
composition corresponds
on the symbol level to the Weyl product, sometimes also called the
twisted product, denoted by $\wpr$. A problem in this field is to find
conditions on the weight functions $\omega _j$ and
$p_j,q_j\in [1,\infty]$, that are necessary and sufficient for the map
\begin{equation}\label{weylmap}
\mathscr S(\rr {2d})\times \mathscr S(\rr {2d}) \ni (a_1,a_2)\mapsto
a_1\wpr a_2 \in \mathscr S(\rr {2d})
\end{equation}
to be uniquely extendable to a map from $M_{(\omega _1)}^{p_1,q_1}(\rr
{2d}) \times M_{(\omega _2)}^{p_2,q_2}(\rr {2d})$ to $M_{(\omega
_0)}^{p_0,q_0}(\rr {2d})$,
which is continuous in the sense that for some constant $C>0$ it holds
\begin{equation}\label{holderyoung3}
\| a_1 \wpr a_2 \|_{M_{(\omega _0)}^{p_0,q_0} }
\leq C \| a_1 \|_{M_{(\omega _1)}^{p_1,q_1} } \| a_2 \|_{M_{(\omega
_2)}^{p_2,q_2} },
\end{equation}
when $a_1 \in M_{(\omega _1)}^{p_1,q_1}(\rr {2d})$ and $a_2 \in
M_{(\omega _2)}^{p_2,q_2}(\rr {2d})$.
Important contributions in this context can be found in
\cite{Gc3,HTW,La1,Sj1,To2}, where Theorem 0.3$'$ in
\cite{HTW} seems to be the most general result so far.

\par

The Weyl product on the Fourier transform side
is given by a twisted convolution, $*_\sigma$. It follows that the
continuity questions here above are the same as finding appropriate
conditions on $\omega _j$ and $p_j,q_j\in [1,\infty]$, in order for
the map
\begin{equation}\label{twistmap}
\mathscr S(\rr {2d})\times \mathscr S(\rr {2d}) \ni (a_1,a_2)\mapsto
a_1*_\sigma a_2 \in \mathscr S(\rr {2d})
\end{equation}
to be uniquely extendable to a map from $W_{(\omega _1)}^{p_1,q_1}(\rr
{2d}) \times W_{(\omega _2)}^{p_2,q_2}(\rr {2d})$ to $W_{(\omega
_0)}^{p_0,q_0}(\rr {2d})$, which is continuous in the sense that for
some constant $C>0$ it holds
\begin{equation}\label{holderyoung4}
\| a_1 *_\sigma a_2 \|_{W_{(\omega _0)}^{p_0,q_0} }
\leq C \| a_1 \|_{W_{(\omega _1)}^{p_1,q_1} } \| a_2 \|_{W_{(\omega
_2)}^{p_2,q_2} },
\end{equation}
when $a_1 \in W_{(\omega _1)}^{p_1,q_1}(\rr {2d})$ and $a_2 \in
W_{(\omega _2)}^{p_2,q_2}(\rr {2d})$. In this context the continuity
result which corresponds to Theorem 0.3$'$ in \cite{HTW} is Theorem
\ref{algthm2} in Section \ref{sec2}.

\par

In the end of Section \ref{sec2} we especially consider the case when
$p_j=q_j=2$. In this case, $W_{(\omega _j)}^{2,2}$ agrees with
$L^2_{(\omega _j)}$, for appropriate choices of $\omega _j$. Hence,
for such $\omega _j$, it
follows immediately from Theorem \ref{algthm2} that the map
\eqref{twistmap} extends to a continuous mapping from $L^2_{(\omega
_1)}(\rr {2d})\times L^2_{(\omega _2)}(\rr {2d})$ to $L^2_{(\omega
_0)}(\rr {2d})$, and that
$$
\nm {a_1*_\sigma a_2}{L^2_{(\omega _0)}}\le C\nm {a_1}{L^2_{(\omega
_1)}}\nm {a_2}{L^2_{(\omega _2)}},
$$
when $a_1\in L^2_{(\omega _1)}(\rr {2d})$ and $a_2\in L^2_{(\omega
_2)}(\rr {2d})$. In Section \ref{sec2} we prove a more general result,
by combining this result with Young's inequality, and then using
interpolation. Finally we use these results in Section \ref{sec3} to
extend the class of possible window functions in the definition of
modulation space norm.

\par

In Sections \ref{sec2}--\ref{sec3}, the second part of the paper, we
cosider spaces of the form $s_{t,p}(\mathscr H_1,\mathscr H_2)$ and
$s_p^A(\mathscr H_1,\mathscr H_2)$, which consist of all
$a\in \mathscr S'(\rr {2d})$ such that $a_t(x,D)$ and $Aa$
respectively are Schatten-von Neumann operator of order $p\in
[1,\infty ]$ from the Hilbert space $\mathscr H_1$ to the Hilbert
space $\mathscr H_2$. From the definitions it follows that the general
continuity results for the usual Schatten-von Neumann classes carry
over to the $s_{t,p}$ and $s_p^A$ spaces, after the usual operator
compositions have been replaced by the Weyl product and twisted
convolution respectively.

\par

In Section \ref{sec4} we consider the case when $\mathscr H_j$ are
equal to $M^2_{(\omega _j)}(\rr d)$ for some appropriate weight
functions $\omega _j$. In this situation we establish Young type
result for dilated convolution and multiplication. More precisely,
assume that
\begin{equation}\label{eq0.4}
{p_1}^{-1} +{p_2}^{-1} = 1 +r^{-1} , \qquad
1 \leq p_1,p_2,r \leq \infty ,
\end{equation}
and that $\mathscr H_j=M^2_{(\omega _j)}$ and $\mathscr
H_{k,j}=M^2_{(\omega _{k,j})}$ are appropriate. If $t_1,t_2\neq 0$ and
$c_1t_1^2+c_2t_2^2=1$, for some choice of $c_1,c_2\in \{ \pm 1\}$,
then we prove that
\begin{alignat*}{3}
a_1(\cdo /t_1)*a_2(\cdo /t_2)&\in s_{t,r}(\mathscr H_1,\mathscr
H_2),&\quad &\text{when}&\quad a_j&\in s_{t,p_j}(\mathscr
H_{1,j},\mathscr H_{2,j}),
\intertext{and}
a_1(t_1\cdo )\cdot a_2(t_2\cdo )&\in s_{t,r}(\mathscr H_1,\mathscr
H_2),&\quad &\text{when}&\quad a_j&\in s_{t,p_j}(\mathscr
H_{1,j},\mathscr H_{2,j}).
\end{alignat*}
In particular, if $\mathscr H_j=\mathscr H_{k,j}=L^2(\rr d)$, then we
recover the results in Section 3 in \cite{To3}. Furthermore, if each
of the operators $a^w_j(x,D)$ are positive with respect to the $L^2$
form, then we prove that the same is true for $(a_1(\cdo
/t_1)*a_2(\cdo /t_2))^w(x,D)$.

\par

%%%%%%%%%%%%%%%%%%%%%%%%%%%%
\section{Preliminaries}\label{sec1}
%%%%%%%%%%%%%%%%%%%%%%%%%%%%

\par

In this section we recall some notations and basic results. The proofs
are in general omitted.

\par

We start by discussing appropriate conditions for the involved weight
functions. Assume that $\omega$ and $v$ are positive and measureable
functions on $\rr d$. Then $\omega$ is called $v$-moderate if
\begin{equation}\label{moderate}
\omega (x+y) \leq C\omega (x)v(y)
\end{equation}
for some constant $C$ which is independent of $x,y\in \rr d$. If $v$
in \eqref{moderate} can be chosen as a polynomial, then $\omega$ is
called polynomially moderated. We let $\mathscr P(\rr d)$ be the set
of all polynomially moderated functions on $\rr d$. If $\omega (x,\xi
)\in \mathscr P(\rr {2d})$ is constant with respect to the
$x$-variable ($\xi$-variable), then we sometimes write $\omega (\xi )$
($\omega (x)$) instead of $\omega (x,\xi )$. In this case we consider
$\omega$ as an element in $\mathscr P(\rr {2d})$ or in $\mathscr P(\rr
d)$ depending on the situation.

\medspace

We recall that the Fourier transform $\mathscr F$ on $\mathscr S'(\rr
d)$, which takes the form \eqref{ftransf} when $f\in L^1(\rr d)$, is a
homeomorphism on $\mathscr S'(\rr d)$ which restricts to a
homeomorphism on $\mathscr S(\rr d)$ and to a unitary operator on
$L^2(\rr d)$.

\par

Let $\fy \in \mathscr S'(\rr d)$ be fixed, and let $f\in \mathscr
S'(\rr d)$. Then the short-time Fourier transform $V_\fy f(x,\xi )$ of
$f$ with respect to the \emph{window function} $\fy$ is the tempered
distribution on $\rr {2d}$ which is defined by
$$
V_\fy f(x,\xi ) \equiv \mathscr F(f \, \overline {\fy (\cdo -x)}(\xi ).
$$
If $f ,\fy \in \mathscr S(\rr d)$, then it follows that
$$
V_\fy f(x,\xi ) = (2\pi )^{-d/2}\int f(y)\overline {\fy
(y-x)}e^{-i\scal y\xi}\, dy .
$$

\medspace

Next we recall some properties on modulation spaces and their Fourier
transforms. Assume that
$\omega \in \mathscr P(\rr {2d})$ and that $p,q\in [1,\infty ]$. Then
the mixed Lebesgue space $L^{p,q}_{1,(\omega )}(\rr {2d})$ consists of
all $F\in L^1_{loc}(\rr {2d})$ such that $\nm F{L^{p,q}_{1,(\omega
)}}<\infty$, and $L^{p,q}_{2,(\omega )}(\rr {2d})$ consists of all
$F\in L^1_{loc}(\rr {2d})$ such that $\nm F{L^{p,q}_{2,(\omega
)}}<\infty$. Here
\begin{align*}
\nm F{L^{p,q}_{1,(\omega )}} &= \Big (\int \Big (\int |F(x,\xi )\omega
(x,\xi )|^p\, dx\Big )^{q/p}\, d\xi \Big )^{1/q},
\intertext{and}
\nm F{L^{p,q}_{2,(\omega )}} &= \Big (\int \Big (\int |F(x,\xi )\omega
(x,\xi )|^q\, d\xi \Big )^{p/q}\, dx \Big )^{1/p},
\end{align*}
with obvious modifications when $p=\infty$ or $q=\infty$.

\par

Assume that $p,q\in [1,\infty ]$, $\omega \in \mathscr P(\rr
{2d})$ and $\fy \in \mathscr S(\rr d)\setminus 0$ are fixed. Then the
\emph{modulation spaces} $M^{p,q}_{(\omega )}(\rr d)$ and
$W^{p,q}_{(\omega )}(\rr d)$ are the Banach spaces which consist of
all $f\in \mathscr S'(\rr d)$ such that
\begin{equation}\label{modnorm}
\nm f{M^{p,q}_{(\omega )}}\equiv \nm {V_\fy f}{L^{p,q}_{1,(\omega
)}}<\infty ,\ \text{and}\ \nm f{W^{p,q}_{(\omega )}}\equiv \nm {V_\fy
f}{L^{p,q}_{2,(\omega )}}<\infty
\end{equation}
respectively. The definitions of $M^{p,q}_{(\omega )}(\rr d)$ and
$W^{p,q}_{(\omega )}(\rr d)$ are independent of the choice of $\fy$
and different $\fy$ gives rise to equivalent norms. (See Proposition
\ref{p1.4} below). From the fact that
$$
V_{\widehat \fy}\widehat f (\xi ,-x) =e^{i\scal x\xi}V_{\check
\fy}f(x,\xi ),\qquad \check \fy (x)=\fy (-x),
$$
it follows that
$$
f\in W^{q,p}_{(\omega )}(\rr d)\quad \Longleftrightarrow \quad
\widehat f\in M^{p,q}_{(\omega _0)}(\rr d),\qquad \omega _0(\xi
,-x)=\omega (x,\xi ).
$$

\par

For conveniency we set $M^p _{(\omega )}= M^{p,p}_{(\omega
)}$, which agrees with $W^p_{(\omega )}=W^{p,p}_{(\omega
)}$. Furthermore we set $M^{p,q}=M^{p,q}_{(\omega )}$ and
$W^{p,q}=W^{p,q}_{(\omega )}$ if $\omega \equiv 1$. If $\omega$ is
given by $\omega (x,\xi )=\omega _1(x)\omega _2(\xi )$, for some
$\omega _1,\omega _2\in \mathscr P(\rr d)$, then $W^{p,q}_{(\omega )}$
is a Wiener amalgam space, introduced by Feichtinger in \cite{Fe2}.

\par

The proof of the following proposition is omitted, since the results
can be found in \cite {Fe3,Fe4,FG1,FG2,FG3,Gc2,To5,
To6, To7, To8}. Here and in what follows, $p'\in[1,\infty]$ denotes the
conjugate exponent of $p\in[1,\infty]$, i.{\,}e. $1/p+1/p'=1$ should
be fulfilled.

\par

\begin{prop}\label{p1.4}
Assume that $p,q,p_j,q_j\in [1,\infty ]$ for $j=1,2$, and $\omega
,\omega _1,\omega _2,v\in \mathscr P(\rr {2d})$ are such that $\omega$
is $v$-moderate and $\omega _2\le C\omega _1$ for some constant
$C>0$. Then the following are true:
\begin{enumerate}
\item[{\rm{(1)}}] if $\fy \in M^1_{(v)}(\rr d)\setminus 0$, then $f\in
M^{p,q}_{(\omega )}(\rr d)$ if and only if \eqref{modnorm} holds,
i.{\,}e. $M^{p,q}_{(\omega )}(\rr d)$. Moreover, $M^{p,q}_{(\omega )}$
is a Banach space under the norm in \eqref{modnorm} and different
choices of $\fy$ give rise to equivalent norms;

\vrum

\item[{\rm{(2)}}] if  $p_1\le p_2$ and $q_1\le q_2$  then
$$
\mathscr S(\rr d)\hookrightarrow M^{p_1,q_1}_{(\omega _1)}(\rr
n)\hookrightarrow M^{p_2,q_2}_{(\omega _2)}(\rr d)\hookrightarrow
\mathscr S'(\rr d)\text ;
$$

\vrum

\item[{\rm{(3)}}] the $L^2$ product $( \cdo ,\cdo )$ on $\mathscr
S$ extends to a continuous map from $M^{p,q}_{(\omega )}(\rr
n)\times M^{p'\! ,q'}_{(1/\omega )}(\rr d)$ to $\mathbf C$. On the
other hand, if $\nmm a = \sup \abp {(a,b)}$, where the supremum is
taken over all $b\in \mathscr {S}(\rr d)$ such that
$\nm b{M^{p',q'}_{(1/\omega )}}\le 1$, then $\nmm {\cdot}$ and $\nm
\cdot {M^{p,q}_{(\omega )}}$ are equivalent norms;

\vrum

\item[{\rm{(4)}}] if $p,q<\infty$, then $\mathscr S(\rr d)$ is dense in
$M^{p,q}_{(\omega )}(\rr d)$ and the dual space of $M^{p,q}_{(\omega
)}(\rr d)$ can be identified
with $M^{p'\! ,q'}_{(1/\omega )}(\rr d)$, through the form $(\cdo  ,\cdo
)_{L^2}$. Moreover, $\mathscr S(\rr d)$ is weakly dense in $M^{\infty
}_{(\omega )}(\rr d)$.
\end{enumerate}
Similar facts hold if the $M^{p,q}_{(\omega )}$ spaces are replaced by
$W^{p,q}_{(\omega )}$ spaces.
\end{prop}

\par

Proposition \ref{p1.4}{\,}(1) allows us  be rather vague concerning
the choice of $\fy \in  M^1_{(v)}\setminus 0$ in
\eqref{modnorm}. For example, if $C>0$ is a
constant and $\mathscr A$ is a subset of $\mathscr S'$, then $\nm
a{M^{p,q}_{(\omega )}}\le C$ for
every $a\in \mathscr A$, means that the inequality holds for some choice
of $\fy \in  M^1_{(v)}\setminus 0$ and every $a\in
\mathscr A$. Evidently, a similar inequality is true for any other choice
of $\fy \in  M^1_{(v)}\setminus 0$, with  a suitable constant, larger
than $C$ if necessary.

\par

In the following remark we list some other properties for modulation
spaces. Here and in what follows we let
$\eabs x =(1+|x|^2)^{1/2}$, when $x\in \rr d$.

\par

\begin{rem}\label{p1.7}
Assume that $p,p_1,p_2,q,q_1,q_2\in [1,\infty ]$ are such that
$$
q_1\le\min (p,p'),\quad q_2\ge \max (p,p'),\quad  p_1\le\min
(q,q'),\quad p_2\ge \max (q,q'),
$$
and that $\omega, v \in
\mathscr P(\rr {2d})$ are such that $\omega$ is $v$-moderate. Then the
following is true:
\begin{enumerate}
\item[{\rm{(1)}}]  if $p\le q$, then $W^{p,q}_{(\omega )}(\rr
d)\subseteq M^{p,q}_{(\omega )}(\rr d)$, and if $p\ge q$, then
$M^{p,q}_{(\omega )}(\rr d)\subseteq W^{p,q}_{(\omega )}(\rr
d)$. Furthermore, if $\omega (x,\xi )=\omega (x)$, then
$$
M^{p,q_1}_{(\omega )}(\rr d)\subseteq W^{p,q_1}_{(\omega )}(\rr d)
\subseteq L^p_{(\omega )}(\rr d)\subseteq W^{p,q_2}_{(\omega )}(\rr
d)\subseteq M^{p,q_2}_{(\omega )}(\rr d).
$$
In particular, $M^2_{(\omega )}=W^2_{(\omega )}=L^2_{(\omega )}$. If
instead $\omega (x,\xi )=\omega (\xi )$, then
$$
W^{p_1,q}_{(\omega )}(\rr d)\subseteq M^{p,q_1}_{(\omega )}(\rr d)
\subseteq \mathscr FL^q_{(\omega )}(\rr d)\subseteq M^{p_2,q}_{(\omega
)}(\rr d)\subseteq W^{p_2,q}_{(\omega )}(\rr d).
$$
Here $\mathscr FL^q_{(\omega _0)}(\rr d)$ consists of all $f\in
\mathscr S'(\rr d)$ such that
$$
\nm {\widehat f \, \omega _0}{L^q}<\infty ;
$$

\vrum

\item[{\rm{(2)}}] if $\omega (x,\xi )=\omega (x)$, then
$$
M^{p,q}_{(\omega )}(\rr d)\subseteq C(\rr d)\quad \Longleftrightarrow
\quad W^{p,q}_{(\omega )}(\rr d)\subseteq C(\rr d)\quad
\Longleftrightarrow \quad q=1.
$$

\vrum

\item[{\rm{(3)}}] $M^{1,\infty}(\rr d)$ and $W^{1,\infty}(\rr d)$ are
convolution algebras. If $C_B'(\rr d)$ is the set of all measures on
$\rr d$ with bounded mass, then
$$
C_B'(\rr d)\subseteq W^{1,\infty}(\rr d)\subseteq M^{1,\infty}(\rr
d)\text ;
$$

\vrum

\item[{\rm{(4)}}] if $x_0\in \rr d$ is fixed and $\omega _0(\xi
)=\omega (x_0,\xi )$, then
$$
M^{p,q}_{(\omega )}\cap \mathscr E' =W^{p,q}_{(\omega )}\cap \mathscr
E' =\mathscr  FL^q_{(\omega _0)}\cap \mathscr E' \text ;
$$

\vrum

\item[{\rm{(5)}}] for each $x,\xi\in \rr d$ and modulation space norm
$\nmm \cdo $ we have
\begin{equation*}
\nmm {e^{i\scal\cdo \xi}f(\cdo -x)}\le C
v(x,\xi) \nmm f,
\end{equation*}
for some constant $C$ which is independent of $f\in \mathscr S'(\rr
d)$;

\vrum

\item[{\rm{(6)}}] if $\tilde {\omega}(x,\xi)=\omega(x,-\xi)$ then
$f\in M^{p,q}_{(\omega)}$ if and only if $\overline f\in
M^{p,q}_{(\tilde {\omega})}$;

\vrum

\item[{\rm{(7)}}] if $s\in \mathbf R$ and $\omega (x,\xi )=\eabs \xi
^s$, then $M^2_{(\omega )}=W^2_{(\omega )}$ agrees with the Sobolev
space $H^2_s$, which consists of all $f\in \mathscr S'$ such that
$\mathscr F^{-1}(\eabs \cdo ^s\widehat f)\in L^2$.
\end{enumerate}
(See e.{\,}g. \cite{Fe3,Fe4,FG1,FG2,FG3,Gc2,To5,
To6,To7,To8}.)
\end{rem}

\par

\begin{rem}\label{usualmod}
Assume that $s,t\in \mathbf R$. In many applications it is common that
functions of the form
$$
\sigma _s(x)=\eabs {x}^s\quad \text{and}\quad \sigma _{s,t}(x,\xi
)\equiv \eabs x^t\eabs \xi ^s,
$$
are involved. Then it easilly follows that $\sigma _s$ and $\sigma
_{s,t}$ are $\sigma _{|s|}$-moderate and $\sigma
_{|s|,|t|}$-moderate respectively. For convenience we set
$$
M^{p,q}_{s,t}=M^{p,q}_{(\sigma _{s,t} )}\quad M^{p,q}_{s}=
M^{p,q}_{(\sigma _s )},
$$
and similarily for other function and distribution
spaces, e.{\,}g. we set $L^p_s = L^p_{(\sigma _s )}$. We note
that for such weight functions we have
$$
M^{p,q}_{s,t}(\rr d)=\sets {f\in \mathscr S'(\rr d)}{\eabs x^t\eabs
D^sf\in M^{p,q}(\rr d)},\quad s,t\in \mathbf R
$$
and
$$
M^{p,q}_s(\rr d) =M^{p,q}_{s,0}(\rr d)\cap M^{p,q}_{0,s}(\rr d),\quad
s\ge 0.
$$
(Cf. \cite {To6}.) In particular, since $M^2=L^2$, it holds
$M^{2}_s=H^2_{s}\cap L^2_{s}$, when $s\ge 0$. (See also \cite{GH1,
Gc2}.)
\end{rem}

\medspace

Next we recall some facts in Chapter XVIII in \cite{Ho1} concerning
pseudo-differential operators. Assume that $a\in \mathscr 
S(\rr {2d})$, and that $t\in \mathbf R$ is fixed. Then the
pseudo-differential operator $a_t(x,D)$ in \eqref{e0.5}
is a linear and continuous operator on $\mathscr S(\rr d)$. For
general $a\in \mathscr S'(\rr {2d})$, the
pseudo-differential operator $a_t(x,D)$ is defined as the continuous
operator from $\mathscr S(\rr d)$ to $\mathscr S'(\rr d)$ with
distribution kernel given by \eqref{atkernel}. This definition makes
sense, since the mappings $\mathscr F_2$ and $F(x,y)\mapsto
F((1-t)x+ty,y-x)$ are homeomorphisms on $\mathscr S'(\rr
{2d})$. Furthermore, by Schwartz kernel theorem it follows that the
map $a\mapsto a_t(x,D)$ is a bijection from $\mathscr S'(\rr {2d})$ to
the set of linear and continuous operators from $\mathscr S(\rr d)$ to
$\mathscr S'(\rr d)$.

\par

We recall that for $s,t\in \mathbf R$ and $a,b\in \mathscr S'(\rr
{2d})$, we have
\begin{equation}\label{calculitransform}
a_{s}(x,D) = b_{t}(x,D) \qquad \Longleftrightarrow \qquad
b(x,\xi )=e^{i(t-s)\scal {D_x }{D_\xi}}a(x,\xi ).
\end{equation}
(Note here that the right-hand side makes sense, since $e^{i(t-s)\scal
{D_\xi }{D_x}}$ on the Fourier transform side is a multiplication by
the bounded function $e^{i(t-s)\scal x \xi }$.

\par

Assume that $t\in \mathbf R$ and $a\in \mathscr S'(\rr {2d})$ are
fixed. Then $a$ is called a rank-one element with respect to $t$, if
the corresponding pseudo-differential operator is of rank-one,
i.{\,}e.
\begin{equation}\label{trankone}
a_t(x,D)f=(f,f_2)f_1, 
\end{equation}
for some $f_1,f_2\in \mathscr S'(\rr d)$. By
straight-forward computations it follows that \eqref{trankone}
is fulfilled, if and only if $a=(2\pi
)^{d/2}W_{f_1,f_2}^{t}$, where  the $W_{f_1,f_2}^{t}$ $t$-Wigner
distribution, defined by the formula
\begin{equation}\label{wignertdef}
W_{f_1,f_2}^{t}(x,\xi ) \equiv \mathscr F(f_1(x+t\cdo
)\overline{f_2(x-(1-t)\cdo )} )(\xi ),
\end{equation}
which takes the form
$$
W_{f_1,f_2}^{t}(x,\xi ) =(2\pi )^{-d/2} \int
f_1(x+ty)\overline{f_2(x-(1-t)y) }e^{-i\scal y\xi}\, dy,
$$
when $f_1,f_2\in \mathscr S(\rr d)$. By combining these facts
with \eqref{calculitransform}, it follows that
\begin{equation}\label{wignertransf}
W_{f_1,f_2}^{t} = e^{i(t-s)\scal {D_x }{D_{\xi}}} W_{f_1,f_2}^{s},
\end{equation}
for each $f_1,f_2\in \mathscr S'(\rr d)$ and $s,t\in \mathbf R$. Since
the Weyl case ($t=1/2$) is particulary important to us, we set
$W_{f_1,f_2}^{t}=W_{f_1,f_2}$ when $t=1/2$. It follows that
$W_{f_1,f_2}$ is the usual (cross-) Wigner distribution of $f_1$ and
$f_2$.

\medspace

Next we recall the definition of symplectic Fourier transform, Weyl
product, twisted convolution and related objects. Assume that $a,b\in
\mathscr S'(\rr {2d})$. Then the Weyl product $a\wpr b$ between $a$
and $b$ is the function or distribution which fulfills $(a\wpr
b)^w(x,D) = a^w(x,D)\circ b^w(x,D)$, provided the right-hand side
makes sense. More general, if $t\in \mathbf R$, then the product $\wpr
_t$ is defined by the formula
\begin{equation}\label{wprtdef}
(a\wpr _t b)_t(x,D) = a_t(x,D)\circ b_t(x,D),
\end{equation}
provided the right-hand side makes sense as a continuous operator from
$\mathscr S(\rr d)$ to $\mathscr S'(\rr d)$.

\par

The even-dimensional vector space $\rr {2d}$ is a (real) symplectic
vector space with the (standard) symplectic form
$$
\sigma(X,Y) = \sigma\big( (x,\xi ); (y,\eta ) \big) = \scal y \xi -
\scal x \eta ,
$$
where $\langle \cdot,\cdot \rangle$ denotes the usual scalar product
on $\rr d$.

\par

The symplectic Fourier transform for $a \in \mathscr{S}(\rr {2d})$
is defined by the formula
\begin{equation*}
(\mathscr{F}_\sigma a) (X)
= \pi^{-d}\int a(Y) e^{2 i \sigma(X,Y)}\,  dY.
\end{equation*}
Then $\mathscr{F}_\sigma^{-1} = {\mathscr{F}_\sigma}$ is continuous on
$\mathscr{S}(\rr {2d})$, and extends as
usual to a homeomorphism on $\mathscr{S}'(\rr {2d})$, and to a
unitary map on $L^2(\rr {2d})$. The symplectic short-time Fourier
transform of $a \in \mathscr{S}'(\rr {2d})$ with respect to the
window function $\fy \in \mathscr{S'}(\rr {2d})$ is defined
by
\begin{equation}\nonumber
\mathcal V_{\fy} a(X,Y) = \mathscr{F}_\sigma \big( a\, \fy (\cdo -X)
\big) (Y),\quad X,Y \in \rr {2d}.
\end{equation}

\par

Assume that $\omega \in \mathscr P(\rr {4d})$. Then we let
$\mathcal M^{p,q}_{(\omega )}(\rr {2d})$ and $\mathcal
W^{p,q}_{(\omega )}(\rr {2d})$ denote the modulation spaces, where
the symplectic short-time Fourier
transform is used instead of the usual short-time Fourier transform in
the definitions of the norms. It follows that any property valid for
$M^{p,q}_{(\omega )}(\rr {2d})$ or $W^{p,q}_{(\omega )}(\rr {2d})$
carry over to $\mathcal M^{p,q}_{(\omega )}(\rr {2d})$ and $\mathcal
W^{p,q}_{(\omega )}(\rr {2d})$ respectively. For example, for the
symplectic short-time Fourier transform we have
\begin{equation}\label{stftsymplfour}
\mathcal V_{\mathscr F_\sigma \fy}(\mathscr F_\sigma a)(X,Y) =
e^{2i\sigma (Y,X)}\mathcal V_\fy a(Y,X),
\end{equation}
which implies that
\begin{equation}\label{twistfourmod}
\mathscr F_\sigma \mathcal M^{p,q}_{(\omega )}(\rr {2d}) = \mathcal
W^{q,p}_{(\omega _0)}(\rr {2d}), \qquad \omega _0(X,Y)=\omega (Y,X).
\end{equation}

\par

Assume that $a,b\in \mathscr S(\rr {2d})$. Then the twisted
convolution of $a$ and $b$ is defined by the formula
\begin{equation}\label{twist1}
(a \ast _\sigma b) (X)
= (2/\pi)^{d/2} \int a(X-Y) b(Y) e^{2 i \sigma(X,Y)}\, dY.
\end{equation}
The definition of $*_\sigma$ extends in different ways. For example,
it extends to a continuous multiplication on $L^p(\rr {2d})$ when $p\in
[1,2]$, and to a continuous map from $\mathscr S'(\rr {2d})\times
\mathscr S(\rr {2d})$ to $\mathscr S'(\rr {2d})$. If $a,b \in
\mathscr{S}'(\rr {2d})$, then $a \wpr b$ makes sense if and only if $a
*_\sigma \widehat b$ makes sense, and then
\begin{equation}\label{tvist1}
a \wpr b = (2\pi)^{-d/2} a \ast_\sigma (\mathscr F_\sigma {b}).
\end{equation}
We also remark that for the twisted convolution we have
\begin{equation}\label{weylfourier1}
\mathscr F_\sigma (a *_\sigma b) = (\mathscr F_\sigma a) *_\sigma b =
\check{a} *_\sigma (\mathscr F_\sigma b),
\end{equation}
where $\check{a}(X)=a(-X)$ (cf. \cite{To1,To3,To4}). A
combination of \eqref{tvist1} and \eqref{weylfourier1} give
\begin{equation}\label{weyltwist2}
\mathscr F_\sigma (a\wpr b) = (2\pi )^{-d/2}(\mathscr F_\sigma
a)*_\sigma (\mathscr F_\sigma b).
\end{equation}

\medspace

Next we consider the operator $A$ in \eqref{Aadef}. We recall that $A$
is a homeomorphism on $\mathscr S'(\rr {2d})$, which restricts to a
homeomorphism on $\mathscr S(\rr
{2n})$ and to unitary map on $L^2(\rr {2d})$. Furthermore, by
Fourier's inversion formula it follows that the inverse is given by
$$
(A^{-1}U)(x,{\xi}) =\mathscr F(U(\cdo /2-x,\cdo /2+x))(-\xi )
$$
when $U\in \mathscr S'(\rr {2d})$, which takes the form of
$$
(A^{-1}U)(x,{\xi}) =(2 {\pi} )^{-d/2}\int e^{i \scal y {\xi}} U(y/2-x,
y/2+x)\, dy ,
$$
when $U\in \mathscr S(\rr {2d})$.

\par

An important reason for considering the operator $Aa$ is its close
connection with the Weyl calculus (cf. Lemma \ref{lemma1.5}), and
that
\begin{equation}\label{Atwistrel}
A(a \ast _{\sigma} b)= Aa \circ Ab,
\end{equation}
when $a\in \mathscr S'(\rr {2d})$ ($a\in \mathscr S(\rr {2d})$) and
$b\in \mathscr S(\rr {2d})$ ($b\in \mathscr S'(\rr {2d})$). (See \cite
{Fo,To1,To3,To4}.) Here we have identified operators with their
kernels.

\par

In the following lemma we list some facts about the operator $A$. The
result is a consequence of Fourier's inversion formula, and the
verifications are left for the reader.

\par

\begin{lemma}\label{lemma1.5}
Let $A$ be as above and let $U=Aa$ where $a \in \mathscr S'(\rr
{2d})$.  Then the following is true:
\begin{enumerate}
\item $\ \check U=A\check a$, if $\check a(X)=a(-X)$;

\vrum

\item $J_{\mathscr F}U = A \mathscr F_\sigma a$, where $J _{\mathscr
F}U(x,y) = U(-x,y)$;

\vrum

\item $A(\mathscr F_\sigma a) = (2\pi )^{d/2}a^w(x,D)$ and
$(a^w(x,D)f,g)=(2\pi )^{-d/2}(Aa,\check g\otimes \overline f)$ when
$f,g\in \mathscr S(\rr d)$;

\vrum

\item the Hilbert space adjoint of $Aa$ equals $A\widetilde a$,
where $\widetilde a(X)=\overline {a(-X)}$. Furthermore, if
$a_1,a_2,b\in \mathscr S(\rr {2d})$, then
\begin{equation*}
(a_1*_\sigma a_2,b) = (a_1,b*_\sigma \widetilde a_2)=(a_2,\widetilde
a_1*_\sigma b),\quad (a_1*_\sigma a_2)*_\sigma b = a_1*_\sigma
(a_2*_\sigma b).
\end{equation*}
\end{enumerate}
\end{lemma}

\par

Next we recall some facts on operators which are positive with respect
to the $L^2$ form (see the end of the introduction). Assume that $T$
is a linear and continuous operator from $\mathscr S(\rr d)$ to
$\mathscr S'(\rr n)$. We say that $T$ is positive semi-definite and
write $T\ge 0$, if $(Tf,f)_{L^2}\ge 0$ for every $f\in \mathscr S(\rr
d)$. Furthermore, we also consider distributions which are positive
with respect to the twisted convolution.

\par

\begin{defn}\label{sigmapositive}
Assume that  $a\in \mathscr D'(\rr {2d})$. Then $a$ is called a
$\sigma$-positive distribution if $(a*_\sigma \fy ,\fy )_{L^2}\ge 0$
for all $\fy \in C_0^\infty (\rr {2d)}$. The set of all
$\sigma$-positive distributions is denoted by $\mathscr S'_+(\rr
{2d})$. Furthermore, $\mathscr S'_+(\rr {2d})\cap C(\rr {2d})$, the
set of $\sigma$-positive (continuous) functions, is denoted by
$C_+(\rr {2d})$.
\end{defn}

\par

In \cite{To4} that it is proved that
$$
\mathscr S'_+\subseteq \mathscr S,\quad \text{and}\quad C_+\subseteq
L^2\cap C_B\cap \mathscr FC_B,
$$
where $C_B(\rr {2d})$ is the set of bounded continuous functions on
$\rr d$ which turns to zero at infinity.

\par

The following result is a straight-forward consequence of the
definitions, and shows that positivity in the sense of Definition
\ref{sigmapositive} is closely related to positive in operator and
pseudo-differential operator theory. (See also \cite{To1,To3}.)

\par

\begin{prop}\label{prop1.11}
Assume that $a\in \mathscr S'(\rr {2d})$. Then
\begin{equation*}
a\in \mathscr S'_+(\rr {2d})\quad \Longleftrightarrow \quad Aa \ge 0
\quad \Longleftrightarrow \quad (\mathscr F_\sigma a)^w(x,D)\ge 0.
\end{equation*}
\end{prop}

\medspace

In the end of Section \ref{sec5} we also consider Toeplitz
operators. Assume that $a\in \mathscr S(\rr {2d})$, $h_1,h_2\in
\mathscr S(\rr d)$, and set $\check f(x)=f(-x)$ when $f$ is a
distribution. Then the Toeplitz operator $\tp _{h_1,h_2}(a)$, with
symbol $a$, and window functions $h_1$ and $h_2$, is defined by the
formula
\begin{equation*}
(\tp _{h_1,h_2}(a)f_1,f_2) = (aV_{\check h_1}f_1,V_{\check h_2}f_2) =
(a(2\cdo )W_{f_1, {h_1}},W_{f_2, {h_2}})
\end{equation*}
when $f_1,f_2\in \mathscr S(\rr d)$. The definition of $\tp
_{h_1,h_2}(a)$ extends in several ways. For example, several
extensions are presented in \cite{CG1,HW,To1,To3,To5,To8,To9}, in such
ways that $h_1,h_2$ and $a$ are permitted to belong to Lebesgue
spaces, modulation spaces or Schatten-von Neumann symbol classes.

\par

The most of these extensions are based on the fact that the
pseudo-differential operator symbol of a Toeplitz opertor can be
expressed by the relation
\begin{equation}\label{ToepWeyl}
\begin{aligned}
(a\ast u)_t(x,D) &= \operatorname{Tp}_{h_1,h_2}(a) \quad \text{with}
\\[1ex]
u(X) &= (2\pi)^{-d/2}W_{h_2,h_1}^t(-X),
\end{aligned}
\end{equation}
when $t=1/2$, $h_1,h_2$ are suitable window functions on $\rr d$ and
$a$ is an appropriate distribution on $\rr {2d}$. (See e.{\,}g. \cite
{CG1,To1,To3,To5,To6,To7,To9}.) For general $t$,
\eqref{ToepWeyl} is an immediate consequence of the case $t=1/2$,
\eqref{wignertransf}, and the fact that
$$
e^{i(t-s)\scal {D_x }{D_{\xi}}} (a\ast u) =a * (e^{i(t-s)\scal {D_x
}{D_{\xi}}} u),
$$
which follows by integrations by parts.

\par

%%%%%%%%%%%%%%%%%%%%%%%%%%%%%%%
\section{Twisted convolution on modulation spaces and Lebesgue
spaces}\label{sec2}
%%%%%%%%%%%%%%%%%%%%%%%%%%%%%%%

\par

In this section we discuss algebraic properties of the twisted
convolution when acting on modulation spaces of the form
$\mathcal W^{p,q}_{(\omega )}$. The most general result is equivalent to
Theorem 0.3$'$ in \cite{HTW}, which concerns continuity for
the Weyl product on modulation spaces of the form $\mathcal
M^{p,q}_{(\omega )}$. Thereafter we use this result to establish
continuity properties for the twisted convolution when acting on
weighted Lebesgue spaces.

\par

The following lemma is important in our investigations. The
proof is omitted since the result is an immediate consequence of Lemma
4.4 in \cite{To2} and its proof, \eqref{stftsymplfour},
\eqref{tvist1} and \eqref{weylfourier1}.

\par

\begin{lemma}\label{cornerstone}
Assume that $a_1\in \mathscr S'(\rr {2d})$, $a_2 \in \mathscr
S(\rr {2d})$ and $\fy _1,\fy _2 \in \mathscr S(\rr {2d})$. Then the
following is true:
\begin{enumerate}
\item if $\fy = \pi^d \fy _1 \wpr \fy _2$, then $\fy \in \mathscr
S(\rr {2d})$, the map
$$
Z\mapsto e^{2 i \sigma(Z,Y)} (\mathcal V_{\chi_1} a_1) (X-Y+Z,Z) \,
(\mathcal V_{\chi_2} a_2)(X+Z,Y-Z)
$$
belongs to $L^1(\rr {2d})$, and
\begin{multline}\label{weylp1}
\mathcal V_{\fy } (a_1 \wpr a_2) (X,Y)
\\[1ex]
=
\int  e^{2 i \sigma(Z,Y)} (\mathcal V_{\chi_1} a_1)(X-Y+Z,Z) \,
(\mathcal V_{\chi_2}a_2)
(X+Z,Y-Z) \, dZ\text ;
\end{multline}

\vrum

\item if $\fy = 2^{-d} \fy _1 *_\sigma \fy _2$, then $\fy \in \mathscr
S(\rr {2d})$, the map
$$
Z\mapsto e^{2 i \sigma(X,Z-Y)} (\mathcal V_{\chi_1} a_1) (X-Y+Z,Z) \,
(\mathcal V_{\chi_2} a_2)(Y-Z,X+Z)
$$
belongs to $L^1(\rr {2d})$, and
\begin{multline}\label{twistp1}
\mathcal V_{\fy } (a_1 *_\sigma a_2) (X,Y)
\\[1ex]
=
\int  e^{2 i \sigma(X,Z-Y)} (\mathcal V_{\chi_1} a_1)(X-Y+Z,Z) \,
(\mathcal V_{\chi_2}a_2)(Y-Z,X+Z) \, dZ
\end{multline}
\end{enumerate}
\end{lemma}

\par

The first part of the latter result is used in \cite{HTW} to prove the
following result, which is essentially a restatement of Theorem 0.3$'$
in \cite{HTW}. Here we assume that the involved weight functions
satisfy
\begin{equation}\label{vikt1}
\omega_0(X,Y) \leq C \omega_1(X-Y+Z,Z)
\omega_2(X+Z,Y-Z),\quad \! \! X,Y,Z\in \rr {2d}.
\end{equation}
for some constant $C>0$, and that $p_j, q_j \in
[1,\infty]$ satisfy 
\begin{align}
\frac {1}{p_1}+\frac {1}{p_2}-\frac {1}{p_0} &= 1-\Big (\frac
{1}{q_1}+\frac {1}{q_2}-\frac {1}{q_0}\Big )\label{pqformulas1}
\intertext{and}
0\le \frac {1}{p_1}+\frac {1}{p_2}-\frac {1}{p_0}&\le \frac
{1}{p_j},\frac {1}{q_j}\le \frac {1}{q_1}+\frac {1}{q_2}-\frac
{1}{q_0} ,\quad j=0,1,2.\label{pqformulas2}
\end{align}

\par

\begin{thm}\label{algthm1}
Assume that $\omega _0,\omega _1,\omega _2\in
\mathscr P(\rr {4d})$ satisfy \eqref{vikt1}, and that $p_j, q_j \in
[1,\infty]$ for $j=0,1,2$, satisfy \eqref{pqformulas1} and
\eqref{pqformulas2}.
Then the map \eqref{weylmap} on  $\mathscr S(\rr{2d})$ extends
uniquely to a continuous map from $\mathcal M_{(\omega
_1)}^{p_1,q_1}(\rr {2d}) \times \mathcal M_{(\omega _2)}^{p_2,q_2}(\rr
{2d})$ to $\mathcal M_{(\omega _0)} ^{p_0,q_0}(\rr {2d})$, and for
some constant $C>0$, the bound \eqref{holderyoung3} holds for every
$a_1\in \mathcal M_{(\omega _1)}^{p_1,q_1}(\rr {2d})$ and $a_2\in
\mathcal M_{(\omega _2)}^{p_2,q_2}(\rr {2d})$.
\end{thm}

\par

The next result is an immediate consequence of  \eqref{twistfourmod},
\eqref{weyltwist2} and Theorem \ref{algthm1}. Here the condition
\eqref{vikt1} should be replaced by
\begin{equation}\label{vikt2}
\omega_0(X,Y) \leq C \omega_1(X-Y+Z,Z)
\omega_2(Y-Z,X+Z),\quad \! \! X,Y,Z\in \rr {2d}.
\end{equation}
and the condition \eqref{pqformulas2} should be replaced by
\begin{equation}\label{pqformulas3}
0\le \frac {1}{q_1}+\frac {1}{q_2}-\frac
{1}{q_0}\le \frac {1}{p_j},\frac {1}{q_j}\le \frac {1}{p_1}+\frac
{1}{p_2}-\frac {1}{p_0},\quad j=0,1,2.
\end{equation}

\par

\begin{thm}\label{algthm2}
Assume that $\omega _0,\omega _1,\omega _2\in
\mathscr P(\rr {4d})$ satisfy \eqref{vikt2}, and that $p_j, q_j \in
[1,\infty]$ for $j=0,1,2$, satisfy \eqref{pqformulas1} and
\eqref{pqformulas3}.
Then the map \eqref{twistmap} on $\mathscr S(\rr {2d})$ extends
uniquely to a continuous map from $\mathcal W_{(\omega
_1)}^{p_1,q_1}(\rr {2d}) \times \mathcal W_{(\omega _2)}^{p_2,q_2}(\rr
{2d})$ to $\mathcal W_{(\omega _0)}^{p_0,q_0}(\rr {2d})$, and for some
constant $C>0$, the bound \eqref{holderyoung4} holds for every $a_1\in
\mathcal W_{(\omega _1)}^{p_1,q_1}(\rr {2d})$ and $a_2\in \mathcal
W_{(\omega _2)}^{p_2,q_2}(\rr {2d})$.
\end{thm}

\par

By using Theorem \ref{algthm2} we may generalize Proposition 1.4 in
\cite{To3} to involve continuity of the twisted convolution on
weighted Lebesgue spaces. Here the condition \eqref{vikt2} is replaced
by
\begin{equation}\label{vikt3}
\omega _0(X_1+X_2) \le C\omega _1(X_1)\omega _2(X_2)
\end{equation}

\par

\begin{thm}\label{twistedleb}
Assume that $\omega _0,\omega _1,\omega _2\in
\mathscr P(\rr {2d})$ and $p,p_1,p_2 \in [1,\infty]$ satisfy
\eqref{vikt3}, $p_1,p_2\le p$ and
\begin{equation*}
\max \Big ( \frac 1p,\frac 1{p'}\Big )  \le \frac
1{p_1}+\frac 1{p_2}-\frac 1{p}\le 1,
\end{equation*}
for some constant $C$. Then the map \eqref{twistmap} extends
uniquely to a continuous mapping from $L^{p_1}_{(\omega _1)}(\rr
{2d})\times L^{p_2}_{(\omega _2)}(\rr {2d})$ to $L^p_{(\omega _0)}(\rr
{2d})$. Furthermore, for some constant $C$ it holds
\begin{multline*}
\nm {a_1*_\sigma a_2}{L^p_{(\omega _0)}} \le C \nm
{a_1}{L^{p_1}_{(\omega _1)}} \nm {a_2}{L^{p_2}_{(\omega _2)}},
\\[1ex]
\text{when}\quad a_1\in L^{p_1}_{(\omega
_1)}(\rr {2d}),\quad \text{and}\quad a_2\in L^{p_2}_{(\omega _2)}(\rr
{2d}).
\end{multline*}
\end{thm}

\par

\begin{proof}
From the assumptions it follows that at most one of $p_1$ and $p_2$
are equal to $\infty$. By reasons of symmetry we may therefore assume
that $p_2<\infty$.

\par

Since $W^2_{(\omega )}=M^2_{(\omega )}=L^2_{(\omega )}$ when $\omega
(X,Y)=\omega (X)$, in view of Theorem 2.2 in \cite{To6}, the result
follows from Theorem \ref{algthm2} in the case $p_1=p_2=p=2$.

\par

Now assume that $1/p_1+1/p_2-1/p=1$, $a_1\in
L^{p_1}(\rr {2d})$ and that $a_2\in \mathscr S(\rr {2d})$. Then
$$
\nm {a_1*_\sigma a_2}{L^p_{(\omega _0)}}\le (2/\pi )^{d/2}\nm {\,
|a_1|* |a_2|\, }{L^p_{(\omega _0)}} \le C \nm {a_1}{L^{p_1}_{(\omega
_1)}}\nm {a_2}{L^{p_2}_{(\omega _2)}},
$$
by Young's inequality. The result now follows in this case from
the fact that $\mathscr S$ is dense in $L^{p_2}_{(\omega _2)}$, when
$p_2<\infty$.

\par

The result now follows in the general case by multi-linear
interpolation between the case $p_1=p_2=p=2$ and the case
$1/p_1+1/p_2-1/p=1$, using Theorem 4.4.1 in \cite{BL} and the fact
that
$$
(L^{p_1}_{(\omega )}(\rr {2d}),(L^{p_2}_{(\omega )}(\rr
{2d}))_{[\theta ]} = L^{p_0}_{(\omega )}(\rr {2d}),\quad
\text{when}\quad \frac
{1-\theta}{p_1}+\frac \theta {p_2} = \frac 1{p_0}.
$$
(Cf. Chapter 5 in \cite{BL}.) The proof is complete.
\end{proof}

\par

By letting $p_1=p$ and $p_2=q\le \min (p,p')$, or $p_2=p$ and
$p_1=q\le \min (p,p')$, Theorem \ref{twistedleb} takes the following
form:

\par

\begin{cor}\label{twistedlebcor}
Assume that $\omega _0,\omega _1,\omega _2\in
\mathscr P(\rr {2d})$ and $p, q \in [1,\infty]$ satisfy \eqref{vikt3},
and $q\le \min (p,p')$ for some constant $C$. Then the map
\eqref{twistmap} extends uniquely to a continuous mapping from
$L^p_{(\omega _1)}(\rr {2d})\times L^q_{(\omega _2)}(\rr {2d})$ or
$L^q_{(\omega _1)}(\rr {2d})\times L^p_{(\omega _2)}(\rr {2d})$ to
$L^p_{(\omega _0)}(\rr {2d})$.
\end{cor}

\par

In the next section we need the following refinement of Theorem
\ref{twistedleb} concerning mixed Lebesgue spaces.

\par

\renewcommand{\rubrik}{Theorem \ref{twistedleb}$'$}

\begin{tom}
Assume that $k\in \{ 1,2\}$, $\omega _0,\omega _1,\omega _2\in
\mathscr P(\rr {2d})$ and $p,p_j, q,q_j \in [1,\infty]$ for $j=1,2$
satisfy \eqref{vikt3}, $p_1,p_2\le p$, $q_1,q_2\le q$ and 
\begin{equation*}
\max \Big ( \frac 1p,\frac 1{p'}, \frac 1q,\frac 1{q'}\Big ) \le
\frac 1{p_1}+\frac 1{p_2}-\frac 1{p}\, ,\, \frac 1{q_1}+\frac 1{q_2}-\frac
1{q}\le 1,
\end{equation*}
for some constant $C$. Then the map \eqref{twistmap} extends
uniquely to a continuous mapping from $L^{p_1,q_1}_{k,(\omega _1)}(\rr
{2d})\times L^{p_2,q_2}_{k,(\omega _2)}(\rr {2d})$ to
$L^{p,q}_{k,(\omega _0)}(\rr {2d})$. Furthermore, for some constant
$C$ it holds
\begin{multline*}
\nm {a_1*_\sigma a_2}{L^{p,q}_{k,(\omega _0)}} \le C \nm
{a_1}{L^{p_1,q_1}_{k, (\omega _1)}} \nm {a_2}{L^{p_2,q_2}_{k,(\omega
_2)}},
\\[1ex]
\text{when}\quad a_1\in L^{p_1,q_1}_{k,(\omega
_1)}(\rr {2d}),\quad \text{and}\quad a_2\in L^{p_2,q_2}_{k,(\omega
_2)}(\rr {2d}).
\end{multline*}
\end{tom}

\par

\begin{proof}
The result follows from Minkowski's inequality when $p_1=q_1=1$ and
when $p_2=q_2=1$. Furthermore, the result follows in the case
$p_1=p_2=q_1=q_2=2$ from Theorem \ref{twistedleb}. In the general
case, the result follows from these cases and multi-linear
interpolation.
\end{proof}

\par

%%%%%%%%%%%%%%%%%%%%%%%%%%%%%%%
\section{Window functions in modulation space norms}\label{sec3}
%%%%%%%%%%%%%%%%%%%%%%%%%%%%%%%

\par

In this section we use the results in the previous section to prove
that the class of permitted windows in the modulation space norms can
be extended. More precisely we have the following.

\par

\begin{prop}\label{possiblewindows}
Assume that $p,p_0,q,q_0\in [1,\infty ]$ and $\omega ,v\in \mathscr
P(\rr {2d})$ are such that $p_0,q_0\le \min (p,p',q,q')$, $\check v=v$
and $\omega$ is $v$-moderate. Also assume that $f\in
\mathscr S'(\rr d)$. Then the following is true:
\begin{enumerate}
\item if $\fy \in M^{p_0,q_0}_{(v)}(\rr d)\setminus 0$, then $f\in
M^{p,q}_{(\omega )}(\rr d)$ if and only if $V_\fy f\in
L^{p,q}_{1,(\omega )}(\rr {2d})$. Furthermore, $\nmm f\equiv \nm
{V_\fy f}{L^{p,q}_{1,(\omega )}}$ defines a norm for $M^{p,q}_{(\omega
)}(\rr d)$, and different choices of $\fy$ give rise to equivalent
norms;

\vrum

\item if $\fy \in W^{p_0,q_0}_{(v)}(\rr d)\setminus 0$, then $f\in
W^{p,q}_{(\omega )}(\rr d)$ if and only if $V_\fy f\in
L^{p,q}_{2,(\omega )}(\rr {2d})$. Furthermore, $\nmm f\equiv \nm
{V_\fy f}{L^{p,q}_{2,(\omega )}}$ defines a norm for $W^{p,q}_{(\omega
)}(\rr d)$, and different choices of $\fy$ give rise to equivalent
norms.
\end{enumerate}
\end{prop} 

\par

For the proof we note that the relation between Wigner distributions
(cf. \eqref{wignertdef} with $t=1/2$) and short-time Fourier is given
by
$$
W_{f,g}(x,\xi )=2^de^{i\scal x\xi /2}V_{\check g}f(2x,2\xi ),
$$
which implies that
\begin{equation}\label{stftwiennorms}
\nm {W_{f,\check \fy}}{L^{p,q}_{k,(\omega _0)}} = 2^d\nm {V_\fy
f}{L^{p,q}_{k,(\omega )}},\quad \text{when}\quad \omega _0(x,\xi
)=\omega (2x,2\xi )
\end{equation}
for $k=1,2$.

\par

Finally, by Fourier's inversion formula it follows that if $f_1,g_2\in
\mathscr S'(\rr d)$ and $f_1,g_2\in L^2(\rr d)$, then
\begin{equation}\label{wigntwconv}
W_{f_1,g_1}*_\sigma W_{f_2,g_2}  = (\check f_2,g_1)_{L^2}W_{f_1,g_2}.
\end{equation}

\par

\begin{proof}[Proof of Theorem \ref{possiblewindows}.]
We may assume that $p_0=q_0=\min (p,p',q,q')$. Assume that $\fy ,\psi
\in M^{p_0,q_0}_{(v)}(\rr d)\subseteq L^2(\rr d)$, where the inclusion
follows from the fact that $p_0,q_0\le 2$ and $v\ge c$ for some
constant $c>0$. Since $\check v= v$, and $\nm {V_\fy \psi }
{L^{p_0,q_0}_{k,(v)}}=\nm {V_\psi \fy}{L^{p_0,q_0}_{k,(v)}}$ when
$\check v=v$, the result follows if we prove that
\begin{equation}\label{eststft3}
\nm {V_\fy f}{L^{p,q}_{k,(\omega )}}\le C\nm {V_\psi
f}{L^{p,q}_{k,(\omega )}}\nm {V_\fy \psi }{L^{p_0,q_0}_{k,(v)}},
\end{equation}
for some constant $C$ which is independent of $f\in \mathscr S'(\rr
d)$ and $\fy ,\psi \in M^{p_0,q_0}_{(v)}(\rr d)$.

\par

If $p_1=p$, $p_2=p_0$, $q_1=q$, $q_2=q_0$, $\omega _0=\omega
(2\cdo )$ and $v_0=v(2\cdo )$, then Theorem
\ref{twistedleb}$'$ and \eqref{wigntwconv} give
\begin{multline*}\label{wigntwconvnorm}
\nm {V_\fy f}{L^{p,q}_{k,(\omega )}} = C_1\nm
{W_{f,\check \fy}}{L^{p,q}_{k,(\omega _0)}}
\\[1ex]
=C_2\nm {W_{f,\check \psi}*_\sigma W_{\psi,\check
\fy}}{L^{p,q}_{k,(\omega _0)}}\le C_3 \nm {W_{f,\check \psi}}
{L^{p,q}_{k,(\omega _0)}}\nm {W_{\psi ,\check
\fy}}{L^{p_0,q_0}_{k,(v_0)}}
\\[1ex]
=C_4 \nm {V_\psi f}{L^{p,q}_{k,(\omega )}}\nm {V_\fy \psi
}{L^{p_0,q_0}_{k,(v)}},
\end{multline*}
and \eqref{eststft3} follows. The proof is complete.
\end{proof}

\par

%%%%%%%%%%%%%%%%%%%%%%%%%%%%%%%%%%%%%%%%%%%%%%%%%%
\section{Schatten-von Neumann classes and pseudo-differential
operators}\label{sec4}
%%%%%%%%%%%%%%%%%%%%%%%%%%%%%%%%%%%%%%%%%%%%%%%%%%

\par

In this section we discuss Schatten-von Neumann classes of
pseudo-differential operators from a Hilbert space $\mathscr H_1$ to
another Hilbert space $\mathscr H_2$. Schatten-von Neumann classes were
introduced by R. Schatten in \cite {Sca} in the case $\mathscr
H_1=\mathscr H_2$. (See also \cite {Si}). The general situation, when
$\mathscr H_1$ is not necessarily equal to $\mathscr H_2$, has
thereafter been considered in \cite{BS,ST,TB}.

\par

Let $\ON (\mathscr H_j)$, $j=1,2$, denote the family of orthonormal
sequences in $\mathscr H_j$, and assume that
$T\, : \, \mathscr H_1\to \mathscr H_2$ is linear, and that $p\in
[1,\infty ]$. Then set
$$
\nm {T}{\mathscr I_p} = \nm {T}{\mathscr I_p(\mathscr H_1,\mathscr
H_2)}\equiv \sup \Big (\sum |(Tf_j,g_j)_{\mathscr H_2}|^p\Big )^{1/p}
$$
(with obvious modifications when $p=\infty$). Here the supremum is
taken over all $(f_j)\in \ON (\mathscr H_1)$ and
$(g_j)\in \ON (\mathscr H_2)$. Then $\mathscr I_p=\mathscr
I_p(\mathscr H_1,\mathscr H_2)$, the Schatten-von Neumann class of
order $p$, consists of all linear and continuous operators $T$ from
$\mathscr H_1$ to $\mathscr H_2$ such that $\nm T{\mathscr
I_p(\mathscr H_1,\mathscr H_2)}$ is finite. We note that $\mathscr
I_\infty (\mathscr H_1,\mathscr H_2)$ agrees with $\mathcal B(\mathscr
H_1,\mathscr H_2)$, the set of linear and continuous operators from
$\mathscr H_1$ to $\mathscr H_2$, with equality in norms. We also let
$\mathcal K(\mathscr H_1,\mathscr H_2)$ be the set of all linear and
compact operators from $\mathscr H_1$ to $\mathscr H_2$, and equip
this space with the operator norm as
usual. (Note that the notation $\mathscr I_\sharp(\mathscr
H_1,\mathscr H_2)$ was used instead of $\mathcal K(\mathscr
H_1,\mathscr H_2)$ in \cite{To8}.) If $\mathscr H_1=\mathscr
H_2$, then the shorter notation $\mathscr I_p(\mathscr H_1)$ is used
instead of $\mathscr I_p(\mathscr H_1,\mathscr H_2)$, and similarily
for $\mathcal B(\mathscr H_1,\mathscr H_2)$ and $\mathcal K(\mathscr
H_1,\mathscr H_2)$.

\par

Assume that $(e_j)$ is an orthonormal basis in $\mathscr H_1$, and
that $S\in \mathscr I_1(\mathscr H_1)$. Then the trace of $S$ is
defined as
$$
\operatorname{tr}_{\mathscr H_1} S =\sum (Se_j,e_j)_{\mathscr H_1}.
$$
For each pairs of operators $T_1,T_2\in \mathscr I_\infty (\mathscr
H_1,\mathscr H_2)$ such that $T_2^*\circ T_1\in \mathscr I_1(\mathscr
H_1)$, the sesqui-linear form
$$
(T_1,T_2)=(T_1,T_2)_{\mathscr H_1,\mathscr H_2} \equiv
\operatorname{tr} _{\mathscr H_1}(T_2^*\circ T_1)
$$
of $T_1$ and $T_2$ is well-defined. Here we note that $T\in \mathscr
I_p(\mathscr H_1,\mathscr H_2)$ if and only if $T^*\in \mathscr
I_p(\mathscr H_2,\mathscr H_1)$. We refer to \cite{Si,TB,BS} for
more facts about Schatten-von Neumann classes.

\medspace

In order for discussing Schatten-von Neumann operators within
the theory of pseudo-differential operators, we assume from now on
that the Hilbert spaces $\mathscr H,\mathscr H_0,\mathscr H_1,\mathscr
H_2,\dots$ are ``tempered'' in the following sense.

\par

\begin{defn}\label{deftempHilb}
The Hilbert space $\mathscr H \subseteq \mathscr S'(\rr d)$ is called
\emph{tempered} (on $\rr d$), if $\mathscr S(\rr d)$ is contained and
dense in $\mathscr H$.
\end{defn}

\par

Assume that $\mathscr H$ is a tempered Hilbert space on $\rr d$. Then
we let $\check {\mathscr H}$ and $\mathscr H^\tau$ be the sets of all
$f\in \mathscr S'(\rr d)$ such that $\check f\in \mathscr H$ and
$\overline f\in \mathscr H$ respectively. Then $\check {\mathscr H}$
and $\mathscr H^\tau$ are tempered Hilbert spaces under the norms
$$
\nm f{\check {\mathscr H}} \equiv \nm {\check f}{\mathscr H}\quad
\text{and}\quad \nm f{\mathscr H^\tau} \equiv \nm {\overline
f}{\mathscr H}
$$
respectively.

\par

The \emph{$L^2$-dual}, $\mathscr H'$, of $\mathscr H$ is the set
of all $\fy \in \mathscr S'(\rr d)$ such that
$$
\nm \fy {\mathscr H'}\equiv \sup |(\fy ,f)_{L^2(\rr d)}|
$$
is finite. Here the supremum is taken over all $f\in \mathscr S(\rr
d)$ such that $\nm f{\mathscr H}\le 1$. Assume that $\fy \in \mathscr
H'$. Since $\mathscr S$ is dense in $\mathscr H$, it follows from the
definitions that the map $f \mapsto (\fy ,f)_{L^2}$ from $\mathscr
S(\rr d)$ to $\mathbf C$ extends uniquely to a continuous mapping from
$\mathscr H$ to $\mathbf C$. The following version of Riesz lemma
is useful for us. In order to be self-contained, we also give a
proof.

\par

\begin{lemma}\label{Hduallemma}
Assume that $\mathscr H \subseteq \mathscr S'(\rr d)$ is a tempered
Hilbert space with $L^2$-dual $\mathscr H'$. Then the following is
true:
\begin{enumerate}
\item $\mathscr H'$ is a tempered Hilbert space which can be
identified with the dual space of $\mathscr H$ through the $L^2$-form;

\vrum

\item there is a unique map $T_{\mathscr H}$ from $\mathscr H$ to
$\mathscr H'$ such that
\begin{equation}\label{defTH}
(f,g)_{\mathscr H}=(T_{\mathscr H}f,g)_{L^2(\rr d)}\text ;
\end{equation}

\vrum

\item if $T_{\mathscr H}$ is the map in {\rm{(2)}}, $(e_j)_{j\in I}$
is an orthonormal basis in $\mathscr H$ and $\ep _j=T_{\mathscr
H}e_j$, then $T_{\mathscr H}$ is isometric, $(\ep _j)_{j\in I}$ is an
orthonormal basis and
$$
(\ep _j,e_k)_{L^2(\rr d)}=\delta _{j,k}.
$$
\end{enumerate}
\end{lemma}

\par

\begin{proof}
We have that $\mathscr S\subseteq \mathscr H'\subseteq \mathscr S'$,
and since $\mathscr S$ is dense in $\mathscr H$, it follows that
$\mathscr S$ is dense also in $\mathscr H'$.

\par

First assume that $f\in \mathscr H$, $g\in \mathscr S(\rr d)$, and
let $T_{\mathscr H}f$ in $\mathscr S'(\rr d)$ be defined by
\eqref{defTH}. By the definitions it follow that $T_{\mathscr H}f\in
\mathscr H'$, and that $T_{\mathscr H}$ from $\mathscr H$ to $\mathscr
H'$ is isometric. Furthermore, since the dual space of $\mathscr H$
can be identified with itself, under the scalar product of $\mathscr H$,
the asserted duality properties of $\mathscr H'$ follow.

\par

Let $(e_j)_{j\in I}$ be an arbitrary orthonormal basis in $\mathscr
H$, and let $\ep _j=T_{\mathscr H}e_j$. Then it follows that $\nm
{\ep _j}{\mathscr H'}=1$ and
$$
(\ep _j,e_k)_{L^2}=(e_j,e_k)_{\mathscr H}=\delta _{j,k}.
$$
Furthermore, if
\begin{alignat*}{2}
f &= \sum \alpha _je_j,&\qquad \fy &= \sum \alpha _j\ep _j
\\[1ex]
g &= \sum \beta _je_j,&\qquad \gamma &= \sum \beta _j\ep _j
\end{alignat*}
are finite sums, and we set $(\fy ,\gamma )_{\mathscr H'}\equiv
(f,g)_{\mathscr H}$, then it follows that $(\cdo ,\cdo )_{\mathscr
H'}$ defines a scalar product on such finite sums in $\mathscr H'$,
and that $\nm \fy{\mathscr H'}^2= (\fy ,\fy )_{\mathscr H'}$. By
continuity extensions it now follows that $(\fy ,\gamma )_{\mathscr
H'}$ extends uniquely to each $\fy ,\gamma \in \mathscr H'$, and that
the identity $\nm \fy{\mathscr H'}^2= (\fy ,\fy )_{\mathscr H'}$
holds. This proves the result.
\end{proof}

\par

In what follows we call the basis $(\ep _j)$ in Lemma \ref{Hduallemma}
as the \emph{dual basis} of $(e_j)$.

\par

\begin{cor}\label{inbHilbmodHilb}
Assume that $\mathscr H$ is a tempered Hilbert space on $\rr d$. Then
$$
M^2_{s,s}(\rr d)\subseteq \mathscr H,\mathscr H' \subseteq
M^2_{-s,-s}(\rr d),
$$
for some $s\ge 0$. Furthermore, $M^2_{s,s}(\rr d)$ is dense in
$\mathscr H$ and $\mathscr H'$, which in turn are dense in
$M^2_{-s,-s}(\rr d)$.
\end{cor}

\par

\begin{proof}
The topology in $\mathscr S$ can be obtained by using the semi-norms
$$
\nm {f}{[s]} \equiv \sum _{|\alpha|,|\beta |\le s}\nm {x^\alpha
D^\beta f}{L^2},\quad s=0,1,2,\dots \, .
$$
From the fact that $\mathscr S$ is continuously embedded in $\mathscr
H$ and in $\mathscr H'$, it therefore follows that
$$
\nm f{\mathscr H}\le C\nm f{[s]} \quad \text{and}\quad \nm \fy{\mathscr
H'}\le C\nm \fy{[s]},
$$
when $f\in \mathscr S$, provided $s$ is chosen large enough.

\par

Since the completion of $\mathscr S(\rr d)$ under $\nm {\cdo}{[s]}$ is
equal to $M^2_{s,s}(\rr d)$, the result follows by a standard argument
of approximation, using the duality propoerties in
Proposition\ref{p1.4} (4), together with  the facts that $\mathscr S$
is dense in $\mathscr H$, $\mathscr H'$, $M^2_{s,s}$ and in
$M^2_{-s,-s}$. The proof is complete.
\end{proof}

\par

Assume that $\mathscr H_1,\mathscr H_2\subseteq \mathscr S'(\rr d)$
are tempered Hilbert spaces, $t\in \mathbf R$ is fixed and that $p\in
[1,\infty]$. Then we let $s_p^A(\mathscr
H_1,\mathscr H_2)$ and $s_{t,p}(\mathscr H_1,\mathscr H_2)$ be the
sets of all $a\in \mathscr S'(\rr {2d})$ such that $Aa\in \mathscr
I_p(\mathscr H_1,\mathscr H_2)$ and $a_t(x,D)\in \mathscr
I_p(\mathscr H_1,\mathscr H_2)$ respectively. We also let $s_\sharp
^A(\mathscr H_1,\mathscr H_2)$ and
$s_{t,\sharp}(\mathscr H_1,\mathscr H_2)$ be the set of all $a\in
\mathscr S'(\rr {2d})$ such that $Aa \in \mathcal K(\mathscr
H_1,\mathscr H_2)$ and $a_t(x,D)\in \mathcal K(\mathscr H_1,\mathscr
H_2)$ respectively. These spaces are equipped by the norms
\begin{alignat*}{3}
\nm a{s_{t,p}(\mathscr H_1,\mathscr H_2)} &\equiv \nm
{a_t(x,D)}{\mathscr I_{p}(\mathscr H_1,\mathscr H_2)},&\qquad \nm
a{s_{p}^A(\mathscr H_1,\mathscr H_2)} &\equiv \nm {Aa}{\mathscr
I_{p}(\mathscr H_1,\mathscr H_2)},
\\[1ex]
\nm a{s_{t,\sharp}(\mathscr H_1,\mathscr H_2)} &\equiv \nm
a{s_{t,\infty}(\mathscr H_1,\mathscr H_2)},&\qquad \nm
a{s_{\sharp}^A(\mathscr H_1,\mathscr H_2)} &\equiv \nm
a{s_{\infty}^A(\mathscr H_1,\mathscr H_2)}.
\end{alignat*}
Since the mappings $a\mapsto Aa$ and $a\mapsto a_t(x,D)$ are
bijections from $\mathscr S'(\rr {2d})$ to the set of linear and
continuous operators from
$\mathscr S(\rr d)$ to $\mathscr S'(\rr d)$, it follows that $a\mapsto
Aa$ and $a\mapsto a_t(x,D)$ restrict to isometric bijections from
$s_p^A(\mathscr H_1,\mathscr H_2)$ and $s_{t,p}(\mathscr H_1,\mathscr
H_2)$ respectively to $\mathscr I_p(\mathscr H_1,\mathscr
H_2)$. Consequently, the properties for $\mathscr I_p(\mathscr
H_1,\mathscr H_2)$ carry over to $s_p^A(\mathscr H_1,\mathscr H_2)$
and $s_{t,p}(\mathscr H_1,\mathscr H_2)$. In particular, elements in
$s_1^A(\mathscr H_1,\mathscr H_2)$ of finite rank (i.{\,}e. elements
of the form $a\in s_1^A(\mathscr H_1,\mathscr H_2)$ such that $Aa$ is
a finite rank operator) are dense in $s_\sharp ^A(\mathscr
H_1,\mathscr H_2)$ and in $s_p^A(\mathscr H_1,\mathscr H_2)$ when
$p<\infty$. Since the Weyl quantization is particularly important in
our considerations we also set
$$
s_{p}^w = s_{t,p}\quad \text{and}\quad s_{\sharp}^w=
s_{t,\sharp},\quad \text{when}\quad t=1/2.
$$

\par

If $\omega _1,\omega _2\in \mathscr P(\rr {2d})$, then we use the
notation $s_p^A(\omega _1,\omega _2)$ instead of $s_p^A(M^2_{(\omega
_1)},M^2_{(\omega _2)})$. Furthermore we set $s_p^A(\omega _1,\omega
_2)=s_p^A(\rr {2d})$ if in addition $\omega _1=\omega _2=1$. In the
same way the notations for $s_{t,p}$, $s_p^w$, $s_{t,\sharp}$
and $s_\sharp ^w$ are simplified.

\par

\begin{rem}\label{calculitransfer}
Assume that $t,t_1,t_2\in \mathbf R$, $p\in [1,\infty ]$, $\mathscr
H_1,\mathscr H_2$ are tempered Hilbert spaces on $\rr d$ and that
$a,b\in \mathscr S'(\rr {2d})$. Then it follows by Fourier's inversion
formula that the map $e^{it\scal {D_x }{D_{\xi}}}$ is a homeomorphism
on $\mathscr S(\rr {2d})$ which extends uniquely to a homeomorphism on
$\mathscr S'(\rr {2d})$. Furthermore, by \eqref{calculitransform} it
follows that $e^{i(t_2-t_1)\scal {D_x
}{D_{\xi}}}$ restricts to an isometric bijection from
$s_{t_1,p}(\mathscr H_1,\mathscr  H_2)$ to $s_{t_2,p}(\mathscr
H_1,\mathscr H_2)$.
\end{rem}

\par

The following proposition shows how $s_{t,p}(\mathscr H_1,\mathscr
H_2)$, $s_p^A(\mathscr H_1,\mathscr H_2)$ and other similar spaces are
linked together. The proof is essentially the same as the proof of
Proposition 5.1 in \cite{TB}. Here and in what follows we let $a^\tau
(x,\xi )=\overline {a(x,-\xi )}$ be the ``torsion'' of $a\in \mathscr
S'(\rr {2d})$.

\par

\begin{prop}\label{identification1}
Assume that $t\in \mathbf R$, $\mathscr H_1,\mathscr H_2$ are tempered
Hilbert spaces in $\rr d$, $a\in \mathscr S'(\rr{2d})$, and that $p\in
[1,\infty ]$. Then
$s_p^w(\mathscr H_1,\mathscr H_2)=s_p^A(\mathscr H_1,\check {\mathscr
H}_2)$. Furthermore, the following conditions are equivalent:
\begin{enumerate}
\item $a\in s_p^w(\mathscr H_1,\mathscr H_2)$;

\vrum

\item $\mathscr F_\sigma a\in s_p^w(\mathscr H_1,\check {\mathscr
H}_2)=s_p^A(\mathscr H_1,\mathscr H_2)$;

\vrum

\item $\overline a\in s_p^w(\mathscr H_2',\mathscr
H_1')$;

\vrum

\item $a^\tau \in s_p^A(\mathscr H_1^\tau ,\mathscr H_2^\tau )$;

\vrum

\item $\check a \in s_p^w(\check {\mathscr H}_1 ,\check {\mathscr H}_2
)$;

\vrum

\item $\widetilde a \in s_p^w(\check {\mathscr H}_2',\check{\mathscr
H}_1')$;

\vrum

\item $e^{i(t-1/2)\scal {D_\xi}{D_x}}a\in s_{t,p}(\mathscr
H_1,\mathscr H_2)$.
\end{enumerate}
\end{prop}

\par

\begin{proof}
Let $a_1=\mathscr F_\sigma a$, $a_2=\overline a$, $a_3=a^\tau$,
$a_4=\check a$ and $a_5=\widetilde a$. Then the equivalences follow
immediately from Remark \ref{calculitransfer} and the equalities
\begin{align*}
(a^w(x,D)f,g) &= (a_1^w(x,D)f,\check g) = (f,a_2^w(x,D)g)
\\[1ex]
&= \overline {(a_3^w(x,D)\overline f,\overline g)} =
(a_4^w(x,D)\check f,\check g) =(\check f,a_5^w(x,D)\check g),
\end{align*}
when $a\in \mathscr S'(\rr {2d})$ and $f,g\in \mathscr S(\rr d)$.
Here the first equality follows from the fact that if $K(x,y)$ is
the distribution kernel of $a^w(x,D)$, then $K(-x,y)$ is the
distribution kernel of $(\mathscr F_\sigma a)^w(x,D)=(2\pi
)^{-d/2}Aa$. (Cf. \cite {To1, To3}.) The proof is complete.
\end{proof}

\par

In Remarks \ref{Schattenrem1} and \ref{Schattenrem2} below we list
some properties which follow from well-known results in the theory of
Schatten-von Neumann classes in combination with \eqref{wprtdef},
\eqref{Atwistrel} and the fact that the mappings $a\mapsto a_t(x,D)$
and $a\mapsto Aa$ are isometric bijections from $s_{t,p}(\mathscr
H_1,\mathscr H_2)$ and $s_{p}^w(\mathscr H_1,\mathscr H_2)$
respectively to $\mathscr I _p(\mathscr H_1,\mathscr H_2)$. (We refer
to \cite {Si,TB,BS} for corresponding results about the $\mathscr
I_p$ spaces.) Here the forms $(\cdo ,\cdo )_{s_{t,2}(\mathscr
H_1,\mathscr H_2)}$ and $(\cdo ,\cdo )_{s_{2}^A(\mathscr H_1,\mathscr
H_2)}$ are defined by the formula
\begin{alignat*}{2}
(a,b)_{s_{t,2}(\mathscr H_1,\mathscr H_2)} &=
(a_t(x,D),b_t(x,D))_{\mathscr I_{2}(\mathscr H_1,\mathscr H_2)},&\quad
a,b &\in s_{t,2}(\mathscr H_1,\mathscr H_2)
\intertext{and}
(a,b)_{s_{2}^A(\mathscr H_1,\mathscr H_2)} &= (Aa,Ab)_{\mathscr
I_{2}(\mathscr H_1,\mathscr H_2)},&\quad a,b &\in s_{2}^A(\mathscr
H_1,\mathscr H_2).
\end{alignat*}
We also recall that $p'\in [1,\infty ]$ is the conjugate exponent for
$p\in [1,\infty ]$, i.{\,}e. $1/p+1/p'=1$. Finally, the set $l^\infty
_0$ consists of all sequences in $l^\infty$ which turns to zero at
infinity, and $l^1_0$ consists of all sequences $(\lambda _j)_{j\in
I}$ such that $\lambda _j=0$ except for finite numbers of $j\in I$.

\par

\begin{rem}\label{Schattenrem1}
Assume that $p,p_j,q,r\in [1,\infty ]$ for $1\le j\le 2$, $t\in
\mathbf R$, and that $\mathscr H_1,\dots ,\mathscr H_4$ are tempered
Hilbert spaces on $\rr d$. Then the following is true:
\begin{enumerate}
\item the set $s_{t,p}(\mathscr H_1,\mathscr H_2)$ is a
Banach space which increases with the parameter $p$. If in addition
$p<\infty$ and $p_1\le p_2$, then $s_{t,p}(\mathscr H_1,\mathscr
H_2)\subseteq s_{t,\sharp}(\mathscr H_1,\mathscr H_2)$,
$s_{t,1}(\mathscr H_1,\mathscr H_2)$ is dense in $s_{t,p}(\mathscr
H_1,\mathscr H_2)$ and in $s_{t,\sharp}(\mathscr H_1,\mathscr H_2)$,
and
\begin{equation}\label{normest}
\nm a{s_{t,p_2}(\mathscr H_1,\mathscr H_2)}\le \nm
a{s_{t,p_1}(\mathscr H_1,\mathscr H_2)}, \quad a\in s_{t,\infty
}(\mathscr H_1,\mathscr H_2)\text ;
\end{equation}

\vrum

\item equality in \eqref{normest} is attained, if and only if $a$ is
of rank one, and then $\nm a{s_{t,p}(\mathscr H_1,\mathscr
H_2)_p}=(2\pi )^{-d/2}\nm {f_0}{\mathscr H_1}\nm {g_0}{\mathscr H_2}$,
when $a$ is given by \eqref{wignertdef};

\vrum

\item if $1/p_1+1/p_2=1/r$, $a_1\in s_{t,p_1}(\mathscr H_1,\mathscr
H_2)$ and $a_1\in s_{t,p_1}(\mathscr H_2,\mathscr H_3)$, then $a_2\wpr
_t a_1\in s_{t,r}(\mathscr H_1,\mathscr H_3)$, and
\begin{equation}\label{holder}
\nm {a_2\wpr _t a_1}{ s_{t,r}(\mathscr H_1,\mathscr H_3)}\le \nm
{a_1}{ s_{t,p_1}(\mathscr H_1,\mathscr H_2)}\nm {a_2}{
s_{t,p_2}(\mathscr H_2,\mathscr H_3)}.
\end{equation}
On the other hand, for any $a\in  s_{t,r}(\mathscr H_1,\mathscr H_3)$,
there are elements $a_1\in  s_{t,p_1}(\mathscr H_1,\mathscr H_2)$ and
$a_2\in  s_{t,p_2}(\mathscr H_2,\mathscr H_3)$ such that $a=a_2\wpr _t
a_1$ and equality holds in \eqref{holder};

\vrum

\item if $\mathscr H_1\subseteq \mathscr H_2$ and $\mathscr
H_3\subseteq \mathscr H_4$, then $s_{t,p}(\mathscr H_2,\mathscr
H_3)\subseteq s_{t,p}(\mathscr H_1,\mathscr H_4)$.
\end{enumerate}
Similar facts hold when the $s_{t,p}$ spaces and the product $\wpr _t$
are replaced by $s_p^A$ spaces and $*_\sigma$.
\end{rem}

\par

\begin{rem}\label{Schattenrem2}
Assume that $p,p_j,q,r\in [1,\infty ]$ for $1\le j\le 2$, $t\in
\mathbf R$, and that $\mathscr H_1,\dots ,\mathscr H_4$ are tempered
Hilbert spaces on $\rr d$. Then the following is true:
\begin{enumerate}
\item the form $(\cdo ,\cdo )_{s_{t,2}(\mathscr H_1,\mathscr H_2)}$ on
$s_{t,1}(\mathscr H_1,\mathscr H_2)$ extends uniquely to a
sesqui-linear and continuous form from $s_{t,p}(\mathscr H_1,\mathscr
H_2)\times s_{t,p'}(\mathscr H_1,\mathscr H_2)$ to $\mathbf C$, and
for every $a_1\in s_{t,p}(\mathscr H_1,\mathscr H_2)$ and $a_2 \in
s_{t,p'}(\mathscr H_1,\mathscr H_2)$, it holds
\begin{align*}
(a_1,a_2)_{s_{t,2}(\mathscr H_1,\mathscr H_2)} &= \overline
{(a_2,a_1)_{s_{t,2}(\mathscr H_1,\mathscr H_2)}},
\\[1ex]
|(a_1,a_2)_{s_{t,2}(\mathscr H_1,\mathscr H_2)}| &\le
\nm {a_1}{s_{t,p}(\mathscr H_1,\mathscr H_2)}\nm
{a_2}{s_{t,p'}(\mathscr H_1,\mathscr H_2)}\quad \text{and}
\\[1ex]
\nm {a_1}{s_{t,p}(\mathscr H_1,\mathscr H_2)} &= \sup |
(a_1,b)_{s_{t,2}(\mathscr H_1,\mathscr H_2)}|,
\end{align*}
where the supremum is taken over all $b \in s_{t,p'}(\mathscr
H_1,\mathscr H_2)$ such that $\nm {b}{s_{t,p'}(\mathscr H_1,\mathscr
H_2)}\le 1$. If in addition $p<\infty$, then the dual space of
$s_{t,p}(\mathscr H_1,\mathscr H_2)$ can be identified with
$s_{t,p'}(\mathscr H_1,\mathscr H_2)$ through this form;

\vrum

\item if $a\in s_{t,\sharp}(\mathscr H_1,\mathscr H_2)$, then
\begin{equation}\label{spektraluppd}
a_t(x,D)f = \sum _{j=1}^\infty \lambda _j(f,f_j)_{\mathscr H_1}g_j,
\end{equation}
holds for some $(f_j)_{j=1}^\infty\in \ON (\mathscr H_1)$,
$(g_j)_{j=1}^\infty\in \ON (\mathscr H_2)$ and $\lambda =(\lambda
_j)_{j=1}^\infty \in l^\infty _0$, where the operator on the
right-hand side of \eqref{spektraluppd} convergences with respect to
the operator norm. Moreover, $a\in
s_{t,p}(\mathscr H_1,\mathscr H_2)$, if and only if $\lambda
\in l^p$, and then
$$
\nm a{s_{t,p}}=(2\pi )^{-d/2}\nm \lambda {l^p}
$$
and the operator on the right-hand side of \eqref{spektraluppd}
converges with respect to the norm $\nm \cdot {s_{t,p}(\mathscr
H_1,\mathscr H_2)}$;

\vrum

\item If $0\le \theta \le 1$ is such that $1/p=(1-\theta )/p_1+\theta
/p_2$, then the (complex) interpolation space
$$
(s_{t,p_1}(\mathscr H_1,\mathscr H_2),s_{t,p_2}(\mathscr H_1,\mathscr
H_2))_{[\theta ]} =s_{t,p}(\mathscr H_1,\mathscr H_2)
$$
with equality in norms.
\end{enumerate}
Similar facts hold when the $s_{t,p}$ spaces are replaced
by $s_p^A$ spaces.
\end{rem}

\par

A problem with the form $(\cdo ,\cdo )_{s_{t,2}(\mathscr H_1,\mathscr
H_2)}$ in Remark \ref{Schattenrem2} is the somewhat complicated
structure. In the following we show that there is a canonical way to
replace this form with $(\cdo ,\co )_{L^2}$. We start with the
following result concerning polar decomposition of compact operators.

\par

\begin{prop}\label{polardecomp}
Assume that $\mathscr H_1$ and $\mathscr H_2$ are tempered Hilbert
spaces on $\rr d$, $a\in s_{t,\sharp}(\mathscr H_1,\mathscr H_2)$ and
that
$p\in [1,\infty ]$. Then
$$
a\equiv \sum _{j\in I}\lambda _jW_{g_j,\fy _j}^t
$$
(with norm convergence) for some orthonormal sequences
$(\fy _j)_{j\in I}$ in $\mathscr H_1'$ and $(f_j)_{j\in I}$ in
$\mathscr H_2$, and a sequence $(\lambda _j)_{j\in I}\in l^\infty _0$
of  non-negative real numbers which decreases to zero at
infinity. Furthermore, $a\in s_{t,p}(\mathscr H_1,\mathscr H_2)$, if
and only if $(\lambda _j)_{j\in I}\in l^p$, and
$$
\nm a{s_{t,p}(\mathscr H_1,\mathscr H_2)} = (2\pi )^{-d/2}\nm
{(\lambda _j)_{j\in I}}{l^p}.
$$
\end{prop}

\par

\begin{proof}
By Remark \ref{Schattenrem2} (2) it follows that if $f\in \mathscr
S(\rr d)$, then
\begin{equation}\label{pseudoexp}
a_t(x,D)f(x)=\sum _{j\in I}\lambda _j(f,f_j)_{\mathscr H_1}g_j
\end{equation}
for some orthonormal sequences $(f_j)$ in $\mathscr H_1$ and $(g_j)$
in $\mathscr H_2$, and a sequence $(\lambda _j)\in l^\infty _0$ of
non-negative real numbers which decreases to zero at infinity. Now let
$(\fy _j)_{j\in I}$ be an orthonormal sequence in $\mathscr H_1'$
such that $(\fy _j,g_k)_{L^2}=\delta _{j,k}$. Then
$(f,f_j)_{\mathscr H_1}=(f,\gamma _j)_{L^2}$, and the result follows
from \eqref{pseudoexp}, and the fact that
$$
(W^t_{g_j,\fy _j})_t(x,D)f= (f,\fy
_j)_{L^2}g_j=(f,f_j)_{\mathscr H_1}g_j.
$$
The proof is complete.
\end{proof}

\par

Next we prove that the duals for $s_{t,p}(\mathscr H_1,\mathscr H_2)$
and $s_{p}^A(\mathscr H_1,\mathscr H_2)$ can be identified with
$s_{t,p'}(\mathscr H_1',\mathscr H_2')$ and $s_{p'}^A(\mathscr
H_1',\mathscr H_2')$ respectively  through the form $(\cdo ,\co
)_{L^2}$.

\par

\begin{thm}\label{schattenidenti}
Assume that $t\in \mathbf R$, $p\in [1,\infty )$ and that $\mathscr H
_1,\mathscr H_2$ are tempered Hilbert spaces on $\rr d$. Then the
$L^2$ form on $\mathscr S(\rr {2d})$ extends
uniquely to a duality between $s_{t,p}(\mathscr H_1,\mathscr H_2)$ and
$s_{t,p'}(\mathscr H_1',\mathscr H_2')$, and the dual space for
$s_{t,p}(\mathscr H_1,\mathscr H_2)$ can be identified with
$s_{t,p'}(\mathscr H_1',\mathscr H_2')$ through this form. Moreover, if
$\ell \in s_{t,p}(\mathscr H_1,\mathscr H_2)^*$ and $a\in
s_{t,p'}(\mathscr H_1',\mathscr H_2')$ are such that $\overline{\ell
(b)}=(a,b)_{L^2}$ when $b\in s_{t,p}(\mathscr H_1,\mathscr H_2)$, then
$$
\nmm \ell = \nm a{s_{t,p'}(\mathscr H_1',\mathscr H_2')}.
$$

\par

The same is true if the $s_{t,p}(\mathscr H_1,\mathscr H_2)$ spaces
are replaced by $s_p^A(\mathscr H_1,\mathscr H_2)$ spaces.
\end{thm}

\par

\begin{proof}
We only prove the assertion in the case $t=1/2$. The general case
follows by similar arguments and is left for the reader. Assume that
$\ell \in s_p^w(\mathscr H_1,\mathscr H_2)^*$. Since the map $b\mapsto
b^w(x,D)$ is an isometric bijection from $s_p^w(\mathscr H_1,\mathscr
H_2)$ to $\mathscr I_p(\mathscr H_1,\mathscr H_2)$, it follows from
Remark \ref{Schattenrem2} (1) that for some $S\in \mathscr
I_{p'}(\mathscr H_1,\mathscr H_2)$ and each orthonormal basis
$(e_j)\in \ON (\mathscr H_1)$ that
\begin{equation}\label{funcident}
\begin{aligned}
\ell (b) &= \operatorname{tr}_{\mathscr H_1}(S^*\circ b^w(x,D))=\sum
(b^w(x,D)e_j,Se_j)_{\mathscr H_2}\quad \text{and}
\\[1ex]
\nmm \ell &=\nm S{\mathscr I_{p'}(\mathscr H_1,\mathscr H_2)},
\end{aligned}
\end{equation}
when $b\in s_p^w(\mathscr H_1,\mathscr H_2)$.

\par

Now let $b\in s_p^w(\mathscr H_1,\mathscr H_2)$ be an arbitrary finite
rank element. Then
$$
b=\sum \lambda _jW_{f_j,\ep _j}\quad \text{and}\quad \nm
b{s_p^w(\mathscr H_1,\mathscr H_2)} = (2\pi )^{-d/2}\nm {(\lambda
_j)}{l^p},
$$
for some orthonormal bases $(\ep _j)\in \ON (\mathscr H_1')$ and
$(f_j)\in \ON (\mathscr H_2)$, and some sequence $(\lambda _j)\in
l^1_0$. We also let $(e_j)\in \ON (\mathscr H_1)$ be the dual basis of
$(\ep _j)$ and $a$ the Weyl symbol of the operator $T_{\mathscr
H_2}\circ S\circ T_{\mathscr H_1'}$. Then $a\in s_{p'}^w(\mathscr
H_1,\mathscr H_2)$ and $\nm a{s_{p'}^w(\mathscr H_1,\mathscr
H_2)}=\nmm \ell$. By straight-forward computations we also get
\begin{multline*}
\ell (b) =\operatorname{tr}_{\mathscr H_1}(S^*\circ b^w(x,D))=\sum
(b^w(x,D)e_j,Se_j)_{\mathscr H_2}
\\[1ex]
=\sum \lambda _j(f_j,Se_j)_{\mathscr H_2} = \sum \lambda
_j(f_j,a^w(x,D)\ep _j)_{L^2(\rr d)}
\\[1ex]
= (2\pi )^{-d/2}\sum \lambda
_j(W_{f_j,\ep _j},a)_{L^2(\rr {2d})}
= (2\pi )^{-d/2}(b,a)_{L^2(\rr {2d})}.
\end{multline*}
Hence $\ell (b) =(2\pi )^{-d/2}(b,a)_{L^2(\rr {2d})}$. The result now
follows from these identities and the fact that the set of finite rank
elements are dense in $s_p^w(\mathscr H_1,\mathscr H_2)$. The proof is
complete.
\end{proof}

\par

Finally we remark that $\mathscr S$ is contained and dense in $s_{t,p}$.

\par

\begin{prop}\label{Schattendense}
Assume that $p\in [1,\infty )$, and that $\mathscr H_1$ and $\mathscr
H_2$ are tempered Hilbert spaces on $\rr d$.  Then $\mathscr S(\rr
{2d})$ is dense in $s_{t,p}(\mathscr H_1,\mathscr H_2)$,
$s_{p}^A(\mathscr H_1,\mathscr H_2)$, $s_{t,\sharp}(\mathscr
H_1,\mathscr H_2)$ and $s_{\sharp}^A(\mathscr H_1,\mathscr
H_2)$. Furthermore,  $\mathscr S(\rr {2d})$ is dense in
$s_{t,\infty}(\mathscr H_1,\mathscr H_2)$ and $s_{\infty}^A(\mathscr
H_1,\mathscr H_2)$ with respect to the weak$^*$ topology.
\end{prop}

\par

\begin{proof}
The result is an immediate consequence of Theorem 4.13 in \cite{To8},
Remark \ref{Schattenrem1} (4) and Corollary \ref{inbHilbmodHilb}. The
proof is complete.
\end{proof}

\par

\begin{rem}
Except for the Hilbert-Schmidt case ($p=2$), it is in general a hard
task to find simple characterizations of Schatten-von Neumann
classes. Important questions therefore concern of finding embeddings
between Schatten-von Neumann classes and well-known function and
distribution spaces. Here we recall some of such embeddings:
\begin{enumerate}
\item[(i)] in Chapter 4 in \cite {Si}, it is proved that if $Q$ is a
unit cube on $\rr d$, $1\le p\le 2$ and $f$ and $g$ are measureable on
$\rr d$ and satisfy
$$
\Big (\sum _{x_\alpha \in \mathbf Z^n} \Big (\int _{x_\alpha +Q}
|f(x)|^2\, dx\Big )^{p/2}\Big )^{1/p}<\infty ,
$$
and similarily for $g$, then $f(x)g(D)\in \mathscr I_p(L^2)$, or
equivalently, $f(x)g(\xi )\in s_{t,p}(\rr {2d})$ when $t=0$;

\vrum

\item[(ii)] Let $B^{p,q}_s(\rr d)$ be the Besov space with parameters
$p,q\in [1,\infty ]$ and $s\in \mathbf R$ (cf. \cite{To3,To5,To7,To8}
for strict definitions). In \cite{To3} sharp embeddings of the form
$$
B_{s_1}^{p,q_1}(\rr {2d})\subseteq s_p^w(\rr {2d})\subseteq
B_{s_2}^{p,q_2}(\rr {2d})
$$
is presented. Here
\begin{equation}\label{q1q2}
q_1=\min (p,p')\quad \text{and}\quad q_2=\max (p,p').
\end{equation}
We also remark that the sharp embedding $B_s^{\infty ,1}(\rr
{2d})\subseteq s_{t,\infty}(\rr {2d})$ for certain choices of $t$ was
proved already in \cite{Bu0,Bu2,M,Su};

\vrum

\item[(iii)] In \cite[Theorem 4.13]{To8} it is proved that if $\omega
_1,\omega _2\in \mathscr P(\rr {2d})$,
$$
\omega (x,\xi ,\eta ,y)=\omega _2(x-ty,\xi +(1-t)\eta )/\omega
_1(x+(1-t)y,\xi -t\eta )
$$
and $p,q_1,q_2\in [1,\infty ]$ satisfy \eqref{q1q2}, then
\begin{equation}\label{inbmodschatt}
M^{p,q_1} _{(\omega )}(\rr {2d})\subseteq s_{t,p}(\omega _1,\omega
_2)\subseteq M^{p,q_2} _{(\omega )}(\rr {2d}).
\end{equation}
In particular, \eqref{inbmodschatt} covers the Schatten-von Neumann
results in \cite{GH1,Sj1,To5}, where similar questions are
considered in the case $\omega _1=\omega _2=\omega =1$. Furthermore,
in \cite {To8}, embeddings between $s_{t,p}(\omega
_1,\omega _2)$ and Besov spaces with $\omega _1=\omega _2$ are
established.
\end{enumerate}
\end{rem}

\par

%%%%%%%%%%%%%%%%%%%%
\section{Young inequalities for weighted Schatten-von Neumann
classes}\label{sec5}
%%%%%%%%%%%%%%%%%%%%

\par

In this section we establish Young type results for dilated
convolutions and multiplications on $s_{t,p}(\mathscr H_1,\mathscr
H_2)$ and on $s_p^A(\mathscr H_1,\mathscr H_2)$, under the assumptions
that $\mathscr H_1$ and $\mathscr H_2$ are appropriate modulation
spaces of Hilbert type.

\par

As a preparation for this we prove some technical lemmas, and start
with the following classification of Hilbert modulation spaces.

\par

\begin{lemma}\label{lemma1.1}
Assume that $\omega \in \mathscr P(\rr {4d})$ is such that $\omega
(x,y,\xi ,\eta )= \omega (x,\xi )$, $\chi \in \mathscr S(\rr d
)\setminus 0$ and that $F\in \mathscr S'(\rr  {2d})$. Then $F\in
M^{2}_{(\omega )}$, if and only if
\begin{equation}\label{form1.1}
\nmm F \equiv \Big (  \iiint |V_\chi (F(\cdo ,y))(x,\xi )\omega (x,\xi
)|^2 \, dxdyd\xi \Big )^{1/2}.
\end{equation}
Furthermore, $F\mapsto \nmm F$ in \eqref{form1.1} defines a norm which
is equivalent to any $M^2_{(\omega )}$ norm.
\end{lemma}

\par

\begin{proof}
We may assume that $\nm \chi {L^2}=1$. Let $\chi _1=\chi \otimes
\chi$, and let $\mathscr F_1F$ denotes the partial Fourier transform
of $F(x,y)$ with respect to the $x$ variable.  By Parseval's formula
we get
\begin{multline*}
\nm F{M^2_{(\omega )}}^2 = \iiiint |(V_{\chi \otimes \chi }F)(x,y,\xi
,\eta )\omega (x,\xi )|^2\, dxdyd\xi d\eta
\\[1ex]
\iint \Big ( \iint |(\mathscr F \big ( F\, \chi _1(\cdo -(x,y))\big
)(\xi ,\eta )\omega (x,\xi )|^2\,   dy d\eta \Big )\, dx d\xi
\\[1ex]
\iint \Big ( \iint |(\mathscr F _1\big (F(\cdo ,z)\, \chi (\cdo
-x)\big )(\xi )\chi (z-y)\omega (x,\xi )|^2\,  dy dz \Big )\, dx d\xi
\\[1ex]
\iint \Big ( \int |(\mathscr F _1\big (F(\cdo ,z)\, \chi (\cdo -x)\big
)(\xi )\omega (x,\xi )|^2\, dz \Big )\, dx d\xi = \nmm F,
\end{multline*}
where the right-hand side is the same as $\nmm F$ in \eqref
{form1.1}. The proof is complete.
\end{proof}

\par

The following lemma is a cornerstone in our further investigations. We
omit the proof since the result agrees with \cite[Lemma 3.2]{To3}.

\par

\begin{lemma}\label{lemma3.2}
Assume that $s,t\in \mathbf R$ satisfies
$(-1)^{j}s^{-2}+(-1)^{k}t^{-2}=1$, for some choice of $j,k\in \{
0,1\}$, and that $a,b\in \mathscr S(\rr {2d})$. Also let $T_{j,z}$ for
$j\in \{ 0,1\}$ and $z\in \rr d$ be the operator on $\mathscr S(\rr
{2d})$, defined by the formula
$$
(T_{0,z}U)(x,y)=(T_{1,z}U)(y,x)=U(x-z,y+z),\quad U\in \mathscr S(\rr
{2d}).
$$
Then
\begin{multline}\label{eq3.3}
A(a(s\cdo )*b(t\cdo ))
\\[1ex]
=(2\pi )^{d/2}\abp {st}^{-d}\int (T_{j,sz}(Aa))(s^{-1}\cdo
)(T_{k,-tz}(Ab))(t^{-1}\cdo )\, dz.
\end{multline}
\end{lemma}

\par

In Theorem \ref{thm3.3} concerns dilated convolutions of $s_{p}^A$
spaces. Here the conditions on the involved weight functions are
\begin{equation}\label{weightcond1}
\begin{aligned}
\vartheta (X_1+X_2) &\le C{\vartheta _{j_1,1}(t_1X_1)}{\vartheta
_{j_2,2}(t_2X_2)}
\\[1ex]
\omega (X_1+X_2) &\le C{\omega _{j_1,1}(t_1X_1)}{\omega
_{j_2,2}(t_2X_2)}
\end{aligned}
\end{equation}
where
\begin{equation}\label{weightcond2}
\omega _{0,k}(X) = \vartheta _{1,k}(-X) = \omega _k(X),\quad \vartheta
_{0,k}(X) = \omega _{1,k}(-X) = \vartheta _k(X)
\end{equation}
and
\begin{equation}\label{eq3.4}
(-1)^{j_1}t_1^{-2} +(-1)^{j_2}t_2^{-2 }=1.
\end{equation}
For convenience we also set $u_t=u(\cdo )$ and $a_{j,t}=a_j(t\cdo )$.

\par

\begin{thm}\label{thm3.3}
Assume that $p_1,p_2,r\in [1,\infty ]$ satisfy 
\eqref{eq0.4}, and that $t_1, t_2\in \mathbf R$ satisfy
\eqref{eq3.4}, for some choices of $j_1,j_2\in \{0,1\}$. Also assume
that $\omega ,\omega _j,\vartheta
,\vartheta _j\in \mathscr P(\rr {2d})$ for $j=1,2$ satisfy
\eqref{weightcond1} and \eqref{weightcond2}. Then the mapping
$(a_1,a_2) \mapsto a_{1,t_1} \ast  a_{2,t_2}$ on $\mathscr S(\rr
{2d})$, extends uniquely to a continuous mapping from
$s^A_{p_1}(1/\omega _1,\vartheta _1) \times s^A_{p_2}(1/\omega
_2,\vartheta _2)$ to $s^A_r(1/\omega ,\vartheta )$. One has the
estimate
\begin{equation}\label{eq3.5}
\nm{ a_{1,t_1} \ast a_{2,t_2}) }{s^A_r(1/\omega
,\vartheta )}\leq C ^d \nm {a_1}{s^A_{p_1}(1/\omega _1,\vartheta
_1)}\nm {a_2}{s^A_{p_2}(1/\omega _2,\vartheta _2)},
\end{equation}
where $C= C_0^2|t_1|^{ -2/p_1}|t_2|^{-2/p_2}$ for some
constant $C_0$ which is independent of $t_1,t_2$ and $d$.
\end{thm}

\par

Before the proof we note that for the involved spaces in Theorem \ref{thm3.3} we have
\begin{equation}\label{compspAspaces}
s^A_p(1/\omega ,\vartheta ) \subseteq s^A_p(\rr {2d})\subseteq s^A_p(\omega ,1/\vartheta ),\quad \text{when}\quad \omega ,\vartheta \ge c,
\end{equation}
for some constant $c>0$. This is an immediate consequence of Remark \ref{Schattenrem1} (4) and that the embeddings $M^{2,2}_{(\omega )}\subseteq M^{2,2}=L^2\subseteq M^{2,2}_{(1/\omega )}$ hold when $\omega$ is bounded from below. In particular,
\begin{equation}\label{compspAspaces2}
s^A_1(1/\omega ,\vartheta ) \subseteq s^A_1(\rr {2d})\subseteq C_B'(\rr {2d})\cap \mathscr FC_B'(\rr {2d})\cap L^2(\rr{2d}),\quad \text{when}\quad \omega ,\vartheta \ge c,
\end{equation}
where the latter embedding follows from Propositions 1.5 and 1.9 in \cite{To4}.

\par

\begin{proof}
Again we only consider the case $j_1=1$
and $j_2=0$, i.{\,}e. $t^{-2}-s^{-2}=1$ when $t_1=s$ and $t_2=t$. The
other cases follows by similar arguments and are left for the
reader. We may assume that $W=T^*\rr d$, and start to prove the
theorem in the case $p=q=r=1$. By Proposition 1.10 and a simple
argument of approximations, it
follows that we may assume that $a_{1}=u$ and $a_{2}=v$ are rank one
elements in $\mathscr S$ and satisfy
$$
\nm u{s^A_1(1/\omega _1,\vartheta _1)}\le C,\qquad \nm
v{s^A_1(1/\omega _2,\vartheta _2)}\le C.
$$
for some constant $C$.
If $A$ is the mapping in (0.1), then it follows that $Au=f_1\otimes
\overline  f_2$ and $Av=g_1\otimes \overline g_2$, and
\begin{align*}
\nm {f_1}{M^2_{(\vartheta _1)}}\nm {f_2}{M^2_{(\omega _1)}}&\le C_1\nm
u{s_1^A(1/\omega _1,\vartheta _1)},
\\[1ex]
\nm {g_1}{M^2_{(\vartheta _2)}}\nm {g_2}{M^2_{(\omega _2)}}&\le C_1\nm
v{s_1^A(1/\omega _2,\vartheta _2)},
\end{align*}
for some vectors $f_1, f_2, g_1, g_2\in \mathscr S$ such that
$$
\nm {f_1}{M^2_{(\vartheta _1)}}\le C_2,\quad \nm {f_2}{M^2_{(\omega
_1)}}\le C_2,\quad \nm {g_1}{M^2_{(\vartheta _2)}}\le C_2,\quad \nm
{g_2}{M^2_{(\omega _2)}}\le C_2,
$$
for some constants $C_1$ and $C_2$.

\par

Set
\begin{equation*}
F(x,z) = f_{2}(x/s+sz) g_{1}(x/t+tz),\qquad
G(y,z) = f_{1}(y/s-sz)g_{2}(y/t-tz).
\end{equation*}
It follows from \eqref{eq3.3} that
\begin{equation*}
A(u_s\ast v_t)(x,y) = (2\pi )^{d/2}|st|^{-d}\int
F(x,z) G(y,z) \, dz.
\end{equation*}
This implies that
\begin{equation}\label{schattennormest1}
\begin{gathered}
\nm{ u_s \ast v_t}{s^A_1(\omega ,\vartheta )} \leq (2 \pi )^{d/2} \abp
{st}^{-d}\int \nm{F(\cdo ,z)}{M^2_{(\vartheta )}}\nm {G(\cdo
,z)}{M^2_{(\omega )}} \, dz
\\[1ex]
\le C\abp {st}^{-d}I_1\cdot I_2,
\end{gathered}
\end{equation}
where
\begin{equation}\label{modnormest3}
\begin{aligned}
I_1 &= \Big ( \iiint |V_\chi (F(\cdo ,z))(x,\xi )\vartheta (x,\xi
)|^2\, dxdzd\xi \Big )^{1/2}
\\[1ex]
I_2 &= \Big ( \iiint |V_\chi (G(\cdo ,z))(x,\xi )\omega  (x,\xi
)|^2\, dxdzd\xi \Big )^{1/2}.
\end{aligned}
\end{equation}
Hence, $I_1\le C\nm F{M^2 _{(\vartheta _0)}}$ and $I_2\le C\nm
G{M^2_{(\omega _0)}}$ by Lemma \ref{lemma1.1}, when $\omega _0(x,y,\xi
,\eta )=\omega (x,\xi )$ and $\vartheta _0(x,y,\xi ,\eta )=\vartheta
(x,\xi )$.

\par

We need to estimate $\nm F{M^2 _{(\vartheta _0)}}$ and $\nm
G{M^2_{(\omega _0)}}$. In order to estimate $\nm F{M^2 _{(\vartheta
_0)}}$ we choose the window function $\chi \in \mathscr S(\rr {2d})$
as
$$
\chi (x,z)=\chi _0(x/s+sz)\chi _0(x/t+tz),
$$
for some real-valued $\chi _0\in \mathscr S(\rr d)$. By taking
$(x_1/s+sz_1,x_1/t+tz_1)$ as new variables when evaluating $V_\chi
F$ we get by formal computations
\begin{multline*}
V_\chi F(x,z,\xi ,\zeta )
\\[1ex]
= (2\pi )^{-d}\iint F(x_1,z_1)\chi
(x_1-x,z_1-z)e^{-i\scal {x_1}{\xi}-i\scal{z_1}{\zeta}}\, dx_1dz_1
\\[1ex]
=(2\pi )^{-d}|st|^{-d}\iint \overline {f_2(x_1)}g_1(z_1)\chi
_0(x_1-(x/s+sz))\chi _0(z_1-(x/t+tz))\times
\\[1ex]
\times e^{-i(\scal
{t^{-1}z_1-s^{-1}x_1}{\xi}+(st)^{-1}\scal{t^{-1}x_1-s^{-1}z_1}{\zeta}}\,
dx_1dz_1
\\[1ex]
=|st|^{-d}\overline {V_{\chi _0}f_2(s^{-1}x+sz,s^{-1}\xi
-(st^2)^{-1}\zeta)} V_{\chi _0}g_1(t^{-1}x+tz,t^{-1}\xi -(s^2t)\zeta
).
\end{multline*}

\par

Furthermore, by \eqref{weightcond1}, \eqref{weightcond2} and the fact that
$t^{-2}-s^{-2}=1$, we obtain
\begin{multline*}
\vartheta (x,\xi ) =\vartheta \big (
(t^{-2}x+z)-(s^{-2}x+z),(t^{-2}\xi -(st)^{-2}\zeta )-(s^{-2}\xi
-(st)^{-2}\zeta )\big )
\\[1ex]
\le C {\omega _1(s^{-1}x+sz,s^{-1}\xi -(st^2)^{-1}\zeta )}
{\vartheta _2(t^{-1}x+tz,t^{-1}\xi -(s^2t)^{-1}\zeta )}
\end{multline*}
A combination of these relations now gives
\begin{equation}\label{stfts}
|V_\chi F(x,z,\xi ,\zeta )\vartheta (x,\xi )| \le C|st|^{-d}J_1\cdot
J_2,
\end{equation}
where
\begin{align*}
J_1 &= |V_{\chi _0}f_2(s^{-1}x+sz,s^{-1}\xi
-(st^2)^{-1}\zeta)\omega _1(s^{-1}x+sz,s^{-1}\xi -(st^2)^{-1}\zeta )|
\intertext{and}
J_2 &= |V_{\chi _0}g_1(t^{-1}x+tz,t^{-1}\xi -(s^2t)\zeta )\vartheta
_2(t^{-1}x+tz,t^{-1}\xi -(s^2t)^{-1}\zeta )|.
\end{align*}

\par

By applying the $L^2$ norm and taking
$$
s^{-1}x+sz,\quad t^{-1}x+tz, \quad s^{-1}\xi -(st^2)^{-1}\zeta ,\quad
t^{-1}\xi -(s^2t)^{-1}\zeta
$$
as new variables of integration we get
\begin{equation}\label{modnormest1}
\nm F{M^2_{(\vartheta )}}\le C|st|^{-2d}\nm {f_2}{M^2_{(\omega
_1)}}\nm {g_1}{M^2_{(\vartheta _2)}}.
\end{equation}

\par

By similar computations it also follows that
\begin{equation}\label{modnormest2}
\nm G{M^2_{(\omega )}}\le C|st|^{-2d}\nm {f_1}{M^2_{(\vartheta
_1)}}\nm {g_2}{M^2_{(\omega _2)}}.
\end{equation}
Hence, a combination of Proposition \ref{polardecomp}, \eqref{schattennormest1}, \eqref{modnormest3},
\eqref{modnormest1} and \eqref{modnormest2} gives
\begin{multline*}
\nm{ u_s \ast v_t}{s^A_1(1/\omega ,\vartheta )}\le C_1|st|^{-d}\nm
{f_1}{M^2_{(\vartheta _1)}}\nm {f_2}{M^2_{(\omega _1)}}\nm
{g_1}{M^2_{(\vartheta _2)}}\nm {g_2}{M^2_{(\omega _2)}}
\\[1ex]
\le C_2|st|^{-d}\nm u{s_1^A(1/\omega _1,\vartheta _1)}\nm v{s_1^A(1/\omega
_2,\vartheta _2)}.
\end{multline*}
This proves the result in the case $p=q=r=1$.

\par

Next we consider the case $p_1=r=\infty$, which implies that $p_2=1$.
Assume that $a\in s^A_\infty (1/\omega _1,\vartheta _1)$ and that
$b,c\in \mathscr S(\rr {2d})$. Then
\begin{equation*}
(a_s*b_t,c)=\abp s^{-4d}(a,{\widetilde b}_{t_0}*c_{s_0}),
\end{equation*}
where $\widetilde b(X)=\overline {b(-X)}$, $s_0=1/s$ and
$t_0=t/s$. We claim that
\begin{equation}\label{bcs1est}
\nm {{\widetilde b}_{t_0}*c_{s_0}}{s_1^A(\omega _1,1/\vartheta _1)}
\le C|s^2/t|^{2d}\nm b{s_1^A(1/\omega _2,\vartheta _2)}\nm
c{s_1^A(\omega ,1/\vartheta )}
\end{equation}

\par

Admitting this for a while, it follows by duality, using Theorem \ref{schattenidenti} that
$$
\nm {a_s*b_t}{s^A_\infty (1/\omega ,\vartheta )}\le
C\abp{s^2/t}^{2d}s^{-4d} \nm a{s^A_\infty (1/\omega _1,\vartheta
_1)} \nm b{s^A_1(1/\omega _2,\vartheta _2)},
$$
which gives \eqref{eq3.5}. The result now follows in the case
$p_1=r=\infty$ and $p_2=1$ from the fact that $\mathscr S$ is dense in
$s^A_1(1/\omega _2,\vartheta _2)$. In the same way the result follows
in the case $p_2=r=\infty$ and $p_1=1$.

\par

For general $p_1,p_2,r\in [1,\infty ]$ the result follows by
multi-linear interpolation, using Theorem 4.4.1 in \cite{BL} and
Remark \ref{Schattenrem2} (3).

\par

It remains to prove \eqref{bcs1est} when $b,c\in \mathscr S(\rr
{2d})$. The condition \eqref{eq3.4} is invariant under the
transformation $(t,s)\mapsto (t_0,s_0)=(t/s,1/s)$. Let
\begin{alignat*}{3}
\widetilde \omega &=1/\omega _1,&\quad \widetilde \vartheta &=
1/\vartheta _1,&\quad \widetilde \omega _1 &=1/\omega ,
\\[1ex]
\widetilde \vartheta _1 &= 1/\vartheta ,&\quad \widetilde \omega _2 &=
\vartheta _2&\quad
\text{and}\quad \widetilde \vartheta _2 &= \omega _2.
\end{alignat*}
If $X_1=-(X+Y)/s$ and $X_2=Y/s$, then it follows that
\begin{align*}
\omega (X_1+X_2) &\le C\vartheta _1(-sX_1)\omega _2(tX_2)
\intertext{and}
\vartheta (X_1+X_2) &\le C\omega _1(-sX_1)\vartheta _2(tX_2),
\intertext{is equivalent to}
\widetilde \omega (X+Y) &\le C\widetilde \vartheta _1(-s_0X)\widetilde
\omega _2(t_0Y)
\intertext{and}
\widetilde \vartheta (X+Y) &\le C\widetilde \omega _1(-s_0X)\widetilde
\vartheta _2(t_0Y).
\end{align*}
Hence, the first part of the proof gives
\begin{multline*}
\nm {\widetilde b_{t_0}*c_{s_0}}{s_1^A(\omega _1,1/\vartheta _1)} =
\nm {\widetilde b_{t_0}*c_{s_0}}{s_1^A(1/\widetilde \omega ,\widetilde
\vartheta )}
\\[1ex]
\le C|s_0t_0|^{-2d}\nm {\widetilde b}{s_1^A(1/\widetilde \omega
_2,\widetilde \vartheta _2)} \nm {\widetilde c}{s_1^A(1/\widetilde
\omega _1,\widetilde \vartheta _1)}
\\[1ex]
= C|s_0t_0|^{-2d}\nm {\widetilde b}{s_1^A(1/\vartheta _2
,\omega _2)} \nm {\widetilde c}{s_1^A(\omega ,1/\vartheta )}
\\[1ex]
= C|s_0t_0|^{-2d}\nm {b}{s_1^A(1/\omega _2
,\vartheta _2)} \nm {\widetilde c}{s_1^A(\omega ,1/\vartheta )},
\end{multline*}
and \eqref{bcs1est} follows.
\end{proof}

\par

\begin{rem}\label{rem3.4}
A proof without any use of interpolation in the case of trivial weight
is presented in Section 2.3 in \cite{To1}.
\end{rem}

\medskip

There is a natural generalization of Theorem \ref{thm3.3} to the case
of more than two factors in the convolution. We recall that the
corresponding Young condition \eqref{eq0.4} for the exponents when we
have convolutions with $N$ functions is
\begin{equation}\tag*{(\ref{eq0.4})$'$}
{p_1}^{-1} +\dots  +{p_N}^{-1} = N-1 +r^{-1} , \qquad
1 \leq p_1,\dots,p_N,r \leq \infty .
\end{equation}
The condition on the involved weight functions is
\begin{equation}\tag*{(\ref{weightcond1})$'$}
\begin{aligned}
\vartheta (X_1+\cdots +X_N) &\le C{\vartheta _{j_1,1}(t_1X_1)}\cdots
{\vartheta _{j_N,N}(t_NX_N)}
\\[1ex]
\omega (X_1+\cdots +X_N) &\le C{\omega _{j_1,1}(t_1X_1)}\cdots {\omega
_{j_N,N}(t_NX_N)}
\end{aligned}
\end{equation}
where
\begin{equation}\tag*{(\ref{weightcond2})$'$}
\omega _{0,k}(X) = \vartheta _{1,k}(-X) = \omega _k(X),\quad \vartheta
_{0,k}(X) = \omega _{1,k}(-X) = \vartheta _k(X)
\end{equation}
and
\begin{equation}\tag*{(\ref{eq3.4})$'$}
(-1)^{j_1}t_1^{-2}+\cdots + (-1)^{j_N}t_N^{-2 }=1.
\end{equation}

\par

\renewcommand{\rubrik}{Theorem \ref{thm3.3}$'$}

\par

\begin{tom}
Assume that $p_1,\dots ,p_N,r\in [1,\infty ]$
satisfy \eqref{eq0.4}$'$, and that $t_1,\dots ,t_N\in \mathbf R$
satisfy \eqref{eq3.4}$'$, for some choices of $j_1,\dots ,j_N\in
\{0,1\}$. Also assume that $\omega ,\omega _j,\vartheta ,\vartheta
_j\in \mathscr P(\rr {2d})$ for $j=1,\dots ,N$ satisfy
\eqref{weightcond1}$'$ and \eqref{weightcond2}$'$.
Then the mapping $(a_1,\dots ,a_N) \mapsto a_{1,t_1} *
\cdots * a_{N,t_N}$ on $\mathscr S(\rr {2d})$ extends uniquely to a
continuous mapping from $s^A_{p_1}(1/\omega _1,\vartheta _1) \times
\cdots \times s^A_{p_N}(1/\omega _N,\vartheta _N)$ to $s^A_r(1/\omega
,\vartheta )$. One has the estimate
\begin{equation}\tag*{(\ref{eq3.5})$'$}
\begin{aligned}
\nm{ a_{1,t_1} \ast \cdots \ast &a_{N,t_N}
}{s^A_r(1/\omega ,\vartheta )}
\\[1ex]
&\leq C ^d\nm {a_1}{s^A_{p_1}(1/\omega _1,\vartheta _1)}\cdots \nm
{a_N}{s^A_{p_N}(1/\omega _N,\vartheta _N)},
\end{aligned}
\end{equation}
where $C = C_0^N|t_1|^{-2/p_1}\cdots |t_N|^{-2/p_N}$ for some
constant $C_0$ which is independent of $N$, $t_1,\dots ,t_N$ and $d$.
\end{tom}

\par

For the proof we need the following lemma.

\par

\begin{lemma}\label{lemmathm3.3'}
Assume that $\rho ,t_1,\dots ,t_N\in \mathbf R\setminus 0$ fulfills
\eqref{eq3.4}$'$ and $\rho ^{-2}+(-1)^{j_N}t_N^{-2}=1$, and set
$t_j'=t_j/\rho$,
\begin{align*}
\widetilde \omega (X) &= \inf \omega _{j_1,1}(t_1'X_1)\cdots  \omega
_{j_{N-1},N-1}(t_{N-1}'X_{N-1})
\intertext{and}
\widetilde \vartheta (X) &= \inf \vartheta _{j_1,1}(t_1'X_1)\cdots
\vartheta _{j_{N-1},N-1}(t_{N-1}'X_{N-1}),
\end{align*}
where the infima are taken over all $X_1,\dots X_{N-1}$ such that
$X=X_1+\cdots X_{N-1}$. Then the following is true:
\begin{enumerate}
\item $\widetilde \omega ,\widetilde \vartheta \in \mathscr P(\rr
{2d})$;

\vrum

\item for each $X_1,\dots X_{N-1}\in \rr {2d}$ it holds
\begin{align*}
\widetilde \omega (X_1+\cdots +X_{N-1}) &= \omega
_{j_1,1}(t_1'X_1)\cdots  \omega _{j_{N-1},N-1}(t_{N-1}'X_{N-1})
\intertext{and}
\widetilde \vartheta (X_1+\cdots +X_{N-1}) &= \vartheta
_{j_1,1}(t_1'X_1)\cdots \vartheta _{j_{N-1},N-1}(t_{N-1}'X_{N-1})
\text ;
\end{align*}

\vrum

\item if $C$ is the same as in \eqref{weightcond1}$'$, then for each
$X,Y\in \rr {2d}$ it holds
$$
\omega (X+Y)\le C\widetilde \omega (\rho X)\omega _N(t_NY),\quad
\text{and}\quad \vartheta (X+Y)\le C\widetilde \vartheta (\rho
X)\vartheta _N(t_NY).
$$
\end{enumerate}
\end{lemma}

\par

\begin{proof}
The assertion (2) follows immediately from the definitions of
$\widetilde \omega$ and $\widetilde \vartheta$, and (3) is an
immediate consequence of \eqref{weightcond1}$'$.

\par

In order to prove (3) we assume that $X=X_1+\cdots +X_{N-1}$.
Since $\omega _{j_1,1}\in \mathscr P(\rr {2d})$, it follows that
\begin{multline*}
\widetilde \omega (X+Y)\le \omega _{j_1,1}(t_1'(X_1+Y))\cdots
\omega _{j_{N-1},N-1}(t_{N-1}'X_{N-1})
\\[1ex]
\le \omega _{j_1,1}(t_1'X_1)\cdots \omega
_{j_{N-1},N-1}(t_{N-1}'X_{N-1})v(Y),
\end{multline*}
for some $v\in \mathscr P(\rr {2d})$. By taking the infimum over all
representations $X=X_1+\cdots +X_N$, the latter inequality becomes
$\widetilde \omega (X+Y)\le \widetilde \omega (X)v(Y)$. This implies
that $\widetilde \omega \in \mathscr P(\rr {2d})$, and in the same way
it follows that $\widetilde \vartheta \in \mathscr P(\rr {2d})$. The
proof is complete.
\end{proof}

\par

\begin{proof}[Proof of Theorem \ref{thm3.3}${}^\prime$.]
We may assume that $N>2$ and that the theorem is already
proved for lower values on $N$.
The condition on $t_j$ is that $c_1t_1^{-2}+\dots +
c_Nt_N^{-2}=1$, where $c_j\in \{ \pm 1\}$. For symmetry reasons we may
assume that $c_1t_1^{-2}+\dots +c_{N-1}t_{N-1}^{-2}={\rho}^{-2}$,
where ${\rho}>0$.
Let $t_j'=t_j/{\rho}$, $\widetilde \omega$ and $\widetilde \vartheta$
be the same as in Lemma \ref{lemmathm3.3'}, and let $r_1\in [1,\infty
]$ be such that $1/r_1+1/p_{N}=1+1/r$. Then
$c_1{(t_1')}^{-2}+\dots +c_{N-1}{(t_{N-1}')}^{-2}=1$, $r_1\ge 1$ since
$p_N\le r$, and
$$
1/p_1+\dots +1/p_{N-1}=N-2+1/r_1.
$$
By the induction hypothesis and Lemma \ref{lemmathm3.3'} (2) it
follows that
$$
b=a_{1,t_1'}*\cdots *a_{N-1,t_{N-1}'} = \rho ^{d(2N-4)}
(a_{1,t_1}*\cdots *a_{N-1,t_{N-1}})(\cdot /\rho )
$$
makes sense as an element in $s_{r_1}^A(1/\widetilde \omega
,\widetilde \vartheta )$, and
\begin{equation*}
\nm {b}{s^A_{r_1}(1/\widetilde \omega ,\widetilde \vartheta )} \le
C\prod _{j=1}^{N-1}\abp {t_j'}^{-2d/p_j}
\nm a{s^A_{p_j}(1/\omega _j,\vartheta _j)},
\end{equation*}
for some constant $C$. Since $1/r_1+1/p_{N}=1+1/r$, it follows from
Lemma \ref{lemmathm3.3'} (3) that $b_\rho * a_{N,t_N}$ makes sense as
an element in $s^A_r(1/\omega ,\vartheta )$, and
\begin{multline*}
\nm {(a_{1,t_1}*\cdots *a_{N-1,t_{N-1}})
*a_{N,t_N}}{s^A_r(1/\omega ,\vartheta )} =
\rho ^{-d(2N-4)}\nm{b_\rho *a_{N,t_N}}{s^A_r(1/\omega ,\vartheta )}
\\[1ex]
\le C_1\nm {a_1}{s^A_{p_1}(1/\omega _1,\vartheta _1)}\cdots \nm
{a_N}{s^A_{p_N}(1/\omega _N,\vartheta _N)},
\end{multline*}
where
\begin{gather*}
C_1=C{\rho}^{d(4-2N-2/r_1)}
\abp {t_N}^{-2d/p_N}\prod _{j=1}^{N-1} \abp
{t_j'}^{-2d/p_j}=C\prod _{j=1}^{N} \abp {t_j}^{-2d/p_j}.
\end{gather*}
This proves the extension assertions. The uniqueness as well as the
symmetri assertions follow from the facts that $\mathscr S$ is dense
in $s_p^A$ when $p<\infty$ and dense in $s_\infty ^A$ with respect to
the weak$^*$ topology, and that at most one $p_j$ is equal to
infinity due to the Young condition. The proof is complete.
\end{proof}

\par

The first part of the following result follows by combining Theorem
\ref{thm3.3}$'$, Proposition \ref{identification1} and
\begin{equation}\label{convfour}
\mathscr F_\sigma (a_1*\cdots *a_N)=\pi ^{dN} (\mathscr F_\sigma
a_1)\cdots (\mathscr F_\sigma a_N),
\end{equation}
when $a_1,\dots ,a_N\in \mathscr S(\rr {2d})$. Here the condition
\eqref{eq3.4}$'$ is replaced by
\begin{equation}\label{eq3.6}
(-1)^{j_1}{t_1^2} +\cdots + (-1)^{j_N}{t_N^2 }=1.
\end{equation}

\par

\begin{thm}\label{thm3.5}
Assume that $p_1,\dots ,p_N,r\in [1,\infty ]$ satisfy
\eqref{eq0.4}$'$, and that $t_1,\dots ,t_N\in \mathbf R$ satisfy
\eqref{eq3.6}, for some choices of $j_1,\dots ,j_N\in \{0,1\}$.
Also assume that $\omega ,\omega _j,\vartheta ,\vartheta _j\in
\mathscr P(\rr {2d})$ for $j=1,\dots ,N$ satisfy
\eqref{weightcond1}$'$ and \eqref{weightcond2}$'$. Then the mapping
$(a_1,\dots ,a_N) \mapsto a_{1,t_1}\cdots a_{N,t_N}$
on $\mathscr S(\rr {2d})$, where $a_{j,t_j}(X)=a_j(t_jX)$, $1\le
j\le N$, extends uniquely to a continuous mapping from
$s^A_{p_1}(1/\omega _1,\vartheta _1) \times \cdots \times
s^A_{p_N}(1/\omega _N,\vartheta _N)$ to $s^A_r(1/\omega
,\vartheta )$. One has the estimate
\begin{equation}\label{eq3.7}
\begin{aligned}
\nm{ a_{1,t_1} \cdots &a_{N,t_N}
}{s^A_r(1/\omega ,\vartheta )}
\\[1ex]
&\leq C^d \nm {a_1}{s^A_{p_1}(1/\omega _1,\vartheta _1)}\cdots \nm
{a_N}{s^A_{p_N}(1/\omega _N,\vartheta _N)}.
\end{aligned}
\end{equation}
where $C= C_0^N|t_1|^{-2/p_1'}\cdots |t_N|^{-2/p_N'}$ for some
constant $C_0$ which is independent of $N$, $t_1,\dots ,t_N$ and $d$.

\par

Moreover, the product is positive semi-definite in the sense of
Definition \ref{sigmapositive}, if this is true for each factor.
\end{thm}

\par

\begin{proof}
When verifying the positivity statement we may argue by induction as
in the proof of Theorem \ref{thm3.3}$'$. This together with
Proposition \ref{prop1.11} and some simple arguments of approximation
shows that it suffices to prove that $a_sb_t$ is positive semi-definite
when $\pm s^2\pm t^2=1$, $st\neq 0$, and $a,b\in \mathscr S(\rr {2d})$
are $\sigma$-positive rank-one element.

\par

We write
\begin{gather*}
a_s {b}_{t} = \pi ^{-d} \mathscr F _\sigma ( \mathscr F _\sigma
{a}_{s} \ast \mathscr F _\sigma {b}_{t})
= \pi ^{-d} \abp {st}^{-2d}\mathscr F _\sigma( (\mathscr F
_\sigma a)_{1/s} \ast (\mathscr F _\sigma b)_{1/t}).
\end{gather*}
If we set for any $U\in \mathscr S (V\oplus V)$,
$$
U_{0,z}(x,y)=U_{1,z}(-y,-x) =U(x+z,y+z),
$$
then it follows from Lemmas \ref{lemma1.5} and
\ref{lemma3.2} that
\begin{equation*}
A({a}_{s} {b}_{t})(x,y)=(2 /{\pi} )^{d/2}\abp {st}^{-d}\int 
(Aa)_{j,z/s}(sx,sy)(Ab)_{k,-z/t}(tx,ty)\, dz,
\end{equation*}
for some choice of $j,k\in \{ 0,1\}$. Since $a,b\in C_+$ are rank-one
elements, it follows that the integrand is of the form ${\phi} _z(x)
\otimes  \overline {{\phi} _z(y)}$ in all these cases. This proves
that $A({a}_{s}{b}_{t})$ is a positive semi-definite operator.
\end{proof}

\par

The following two theorems follow immediately from Proposition
\ref{identification1}, Theorem \ref{thm3.3}$'$ and Theorem
\ref{thm3.5}. Here the condition \eqref{weightcond2}$'$ is replaced by
\begin{equation}\tag*{(\ref{weightcond2})$''$}
\omega _{0,k}(X) = \vartheta _{1,k}(X) = \omega _k(X),\quad \vartheta
_{0,k}(X) = \omega _{1,k}(X) = \vartheta _k(X).
\end{equation}

\par

\begin{thm}\label{thm3.3weyl}
Assume that $p_1,\dots ,p_N,r\in [1,\infty]$
satisfy \eqref{eq0.4}$'$, and that $t_1,\dots ,t_N\in \mathbf R$
satisfy \eqref{eq3.4}$'$, for some choices of $j_1,\dots ,j_N\in
\{0,1\}$. Also assume that $\omega ,\omega _j,\vartheta ,\vartheta
_j\in \mathscr P(\rr {2d})$ for $j=1,\dots ,N$ satisfy
\eqref{weightcond1}$'$ and \eqref{weightcond2}$''$.
Then the mapping $(a_1,\dots ,a_N) \mapsto a_{1,t_1} *
\cdots * a_{N,t_N}$ on $\mathscr S(\rr {2d})$, where
$a_{j,t_j}(X)=a_j(t_jX)$, $1\le j\le N$, extends uniquely to a
continuous mapping from
$$
s_{p_1}^w(1/\omega _1,\vartheta _1) \times
\cdots \times s_{p_N}^w(1/\omega _N,\vartheta _N)
$$
to $s_{r}^w(1/\omega ,\vartheta )$. One has the estimate
\begin{equation}\label{eq3.5pseudo}
\begin{aligned}
\nm{ a_{1,t_1} \ast \cdots \ast &a_{N,t_N}
}{s^w_r(1/\omega ,\vartheta )}
\\[1ex]
&\leq C ^d\nm {a_1}{s_{p_1}^w(1/\omega _1,\vartheta _1)}\cdots \nm
{a_N}{s_{p_N}^w(1/\omega _N,\vartheta _N)},
\end{aligned}
\end{equation}
where $C = C_0^N|t_1|^{-2/p_1}\cdots |t_N|^{-2/p_N}$ for some
constant $C_0$ which is independent of $N$, $t_1,\dots ,t_N$ and $d$.

\par

Moreover, if $a_j^w(x,D)\ge 0$ for each $1\le j\le N$, then
$(a_{1,t_1} \ast \cdots \ast a_{N,t_N})^w(x,D)\ge 0$.
\end{thm}

\par

\begin{thm}\label{thm3.5weyl}
Assume that $p_1,\dots ,p_N,r\in [1,\infty ]$ satisfy
\eqref{eq0.4}$'$, and that $t_1,\dots ,t_N\in \mathbf R$ satisfy
\eqref{eq3.6}, for some choices of $j_1,\dots ,j_N\in \{0,1\}$. Also
assume that $\omega ,\omega _j,\vartheta ,\vartheta _j\in \mathscr
P(\rr {2d})$ for $j=1,\dots ,N$ satisfy \eqref{weightcond1}$'$ and
\eqref{weightcond2}$''$. Then the mapping $(a_1,\dots ,a_N) \mapsto
a_{1,t_1} \cdots a_{N,t_N}$ on $\mathscr S(\rr {2d})$, where
$a_{j,t_j}(X)=a_j(t_jX)$, $1\le j\le N$, extends uniquely to a
continuous mapping from
$$
s_{p_1}^w(1/\omega _1,\vartheta _1) \times
\cdots \times s_{p_N}^w(1/\omega _N,\vartheta _N)
$$
to $s_{r}^w(1/\omega ,\vartheta )$. One has the estimate
\begin{equation}\label{eq3.5weyl}
\begin{aligned}
\nm{ a_{1,t_1} \cdots &a_{N,t_N}
}{s^w_r(1/\omega ,\vartheta )}
\\[1ex]
&\leq C^d \nm {a_1}{s_{p_1}^w(1/\omega _1,\vartheta _1)}\cdots \nm
{a_N}{s_{p_N}^w(1/\omega _N,\vartheta _N)},
\end{aligned}
\end{equation}
where $C = C_0^N|t_1|^{-2/p_1'}\cdots |t_N|^{-2/p_N'}$ for some
constant $C_0$ which is independent of $N$, $t_1,\dots ,t_N$ and $d$.
\end{thm}

\par

\begin{rem}\label{thm3.3pseudo}
Theorem \ref{thm3.3weyl} can also be generalized to involve $s_{t,p}$
spaces, for general $t\in \mathbf R$.

\par

In fact, assume that $p_j$, $r$, $t_j$, $\omega$, $\omega _j$,
$\vartheta$ and $\vartheta _j$ for $1\le j\le N$ are the same as in
Theorems \ref{thm3.3weyl} and \ref{thm3.5weyl}. Also assume that
$t\in \mathbf R$, and let $\tau _k=t$ when $j_k=0$ and $\tau _k=1-t$
when $j_k=1$. (The numbers $j_k$ are the same as in \eqref{eq3.4}$'$.)

\par

Then the mapping $(a_1,\dots ,a_N) \mapsto a_{1,t_1}
\ast \cdots a_{N,t_N}$ on $\mathscr S(\rr {2d})$, extends uniquely to a
continuous mapping from
$$
s_{\tau _{1},p_1}(1/\omega _1,\vartheta _1) \times
\cdots \times s_{\tau _{N},p_N}(1/\omega _N,\vartheta _N)
$$
to $s_{t,r}(1/\omega ,\vartheta )$. Furthermore it holds
\begin{equation}\label{eq3.5pseudo2}
\begin{aligned}
\nm{ a_{1,t_1} \ast \cdots \ast &a_{N,t_N}
}{s^A_r(1/\omega ,\vartheta )}
\\[1ex]
&\leq C^d \nm {a_1}{s_{\tau _1,p_1}(1/\omega _1,\vartheta _1)}\cdots \nm
{a_N}{s_{\tau _N,p_N}(1/\omega _N,\vartheta _N)}.
\end{aligned}
\end{equation}
where $C = C_0^N|t_1|^{-2a/p_1}\cdots |t_N|^{-2/p_N}$ for some
constant $C_0$ which is independent of $N$, $t_1,\dots ,t_N$ and $d$.

\par

Moreover, if $(a_j)_{\tau _j}(x,D)\ge 0$ for each $1\le j\le N$, then
$(a_{1,t_1} \ast \cdots \ast a_{N,t_N})_t(x,D)\ge 0$.

\medspace

When proving this we first assume that $a_1,\dots ,a_N\in
\mathscr S$. By Proposition \ref{identification1} we get
\begin{multline*}
\nm{ a_{1,t_1} \ast \cdots \ast a_{N,t_N}}{s_{t,r}(1/\omega ,\vartheta
)} = \nm{ e^{-i(t-1/2)\scal {D_x}{D_\xi}}(a_{1,t_1} \ast \cdots \ast
a_{N,t_N})}{s^w_r(1/\omega ,\vartheta )}
\\[1ex]
=\nm {b_1*\cdots *b_N}{s^w_r(1/\omega ,\vartheta )},
\end{multline*}
where
$$
b_k = e^{-i(-1)^{j_k}(t-1/2)\scal {D_x}{D_\xi}/t_k^2}(a_k(t_k\cdo )) =
(e^{-i(-1)^{j_k}(t-1/2)\scal {D_x}{D_\xi}}a_k)(t_k\cdo ).
$$
Hence by Theorem \ref{thm3.3weyl} we get
$$
\nm{ a_{1,t_1} \ast \cdots \ast a_{N,t_N}}{s_{t,r}(1/\omega ,\vartheta
)} \le CI_1\cdots I_N,
$$
where
\begin{equation*}
I_k = \nm {e^{-i(-1)^{j_k}(t-1/2)\scal
{D_x}{D_\xi}}a_k}{s^w_{p_k}(1/\omega _k,\vartheta _k)}
= \nm {a_k}{s_{\tau _k,p_k}(1/\omega _k,\vartheta _k)}.
\end{equation*}
This gives \eqref{eq3.5pseudo2}.

\par

The result now follows from \eqref{eq3.5pseudo2} and the fact that
$\mathscr S$ is dense in $s_{t,p}(\omega _1,\omega _2)$ when
$p<\infty$, and dense in $s_{t,\infty}(\omega _1,\omega _2)$ with
respect to the weak$^*$ topology.
\end{rem}

\par

Next we consider elements in $s^A_1(1/v ,v )$, where $v =\check
v\in \mathscr P(\rr {2d})$ is submultiplicative. We note that each element
in $s^A_1(1/v ,v )$ is a continuous function which turns to zero at infinity,
since \eqref{compspAspaces2} shows that $s^A_1(1/v ,v )\subseteq C_B(\rr {2d})$.

\par

It follows that any product of odd numbers of elements in
$s^A_1(1/v ,v )$ are again in $s^A_1(1/v ,v)$. In fact, assume that
$a_1,\dots ,a_N\in s^A_1(1/v ,v)$, $|\alpha |$ is odd, and that $t_j=1$.
Then it follows from Theorem \ref{thm3.5} that  $a_1^{\alpha _1}\cdots
a_N^{\alpha _N}\in s^A_1(1/v ,v )$, and
\begin{equation}\label{monomschatten}
\nm {a_1^{\alpha _1}\cdots
a_N^{\alpha _N}}{s^A_1(1/v ,v )}\le C_0^{d|\alpha |} \prod
\nm {a_j}{s^A_1(1/v ,v )}^{\alpha _j},
\end{equation}
for some constant $C_0$ which is independent of $\alpha$ and $d$.

\par

Furthermore, if in addition $a_1,\dots ,a_N$ are $\sigma$-positive in
the sense of Definition \ref{sigmapositive}, then the same is true for
$a_1^{\alpha _1}\cdots a_N^{\alpha _N}$. The following result is an
immediate consequence of these observations.

\par

\begin{prop}\label{prop3.10}
Assume that $a_1,\dots ,a_N\in s^A_1(1/v ,v )$, where
$v=\check v \in \mathscr P(\rr {2d})$ is submultiplicative,
$C_0$ is the same as in \eqref{monomschatten}, and assume that
$R_1,\dots ,R_N>0$. Also assume that $f,g$ are odd analytic functions
from the polydisc
$$
\sets {z\in \mathbf C^N}{\abp {z_j}<C_0R_j}
$$
to $\mathbf C$, with expansions
$$
f(z)=\sum _\alpha c_\alpha z^\alpha \quad \text{and}\quad g(z)=\sum
_\alpha \abp {c_\alpha } z^\alpha .
$$
Then $f(a)=f(a_1,\dots ,a_N)$ is
well-defined and belongs to $s^A_1(1/v ,v )$. One has the estimate
\begin{equation*}
\nm {f(a)}{s^A_1(1/v ,v )}\le g( C_0\nm
{a_1}{s^A_1(1/v ,v )},\dots ,C_0\nm {a_N}{s^A_1(1/v,v)} ).
\end{equation*}
If in addition $a_1,\dots ,a_N\in C_+(\rr {2d})$, then $g(a)\in C_+(\rr {2d})$.
\end{prop}

\par

An open question for the author is whether the Theorems
\ref{thm3.3}--\ref{thm3.5weyl} and Remark \ref{thm3.3pseudo} are true
for other dilations. This might then lead
to improvements of Proposition \ref{prop3.10}. In this context we note
that $s^A_1(\rr {2d})$, and therefore $s^A_\infty (\rr {2d})$ by
duality, are not stable under dilations (see Proposition 2.1.12 in
\cite {To1} or Proposition 5.4 in \cite{To4}). We refer to
\cite{To3,To4} for a further properties of $\sigma$-positive
functions and distributions.

For rank one elements we also have the following positivity result.

\par

\begin{prop}\label{prop4.10}
Assume that $v ,v_1\in \mathscr P(\rr {2d})$ are
submultiplicative and fulfill $v_1=v(\cdo /\sqrt 2\, )$,
$u\in s^w_\infty (1/\omega ,\omega )$ is an
element of rank one, and let $a(X)=\abp {u(X/\sqrt 2)}^2$. Then $a\in
s^w_1 (1/v_1,v_1)$, and $a^w(x,D)\ge 0$.
\end{prop}

\par

\begin{proof}
Since $u$ is rank one, it follows from Proposition
\ref{identification1} that $u,\overline u\in s^w_1 (1/v
,v)$, which implies that $a\in s^w_1 (1/v_1
,v_1)$ in view of Theorem \ref{thm3.5weyl}. The result now follows from
this fact and Proposition 4.10 in \cite{To3}.
\end{proof}

\medspace

We finish the section by applying our results on Toeplitz
operators. The following result is can be considered as a parallel
result to the recent results in \cite{TB}, especially to Theorem 3.1
and Theorem 3.5 in \cite{TB}. It also generalizes Proposition 4.5 in
\cite{To5}.

\par

\begin{thm}\label{cont3}
Assume that $p\in [1,\infty ]$ and $\omega ,\omega _0,\vartheta
,\vartheta _j\in \mathscr P(\rr {2d})$ for $j=0,1,2$ satisfy
\begin{align*}
\omega (X_1-X_2) &\le C\omega _0(\sqrt 2\, X_1)\vartheta _2(X_2)
\intertext{and}
\vartheta (X_1-X_2) &\le C\vartheta _0(\sqrt 2\, X_1)\vartheta _1(X_2)
\end{align*}
Then the definition of $\tp _{h_1,h_2}(a)$ extends uniquely to each
$a\in \mathscr S'(\rr {2d})$ and $h_j\in M^2_{(\vartheta _j)}$ for $j=1,2$
such that $a(\sqrt 2\cdo )\in s^w_p(1/\omega _0,\vartheta _0)$, and for
some constant $C$ it holds
$$
\nm {\tp _{h_1,h_2}(a)}{\mathscr I_p(M^2_{(1/\omega )},M^2_{(\vartheta
)})}\le C\nm {a(\sqrt 2\cdo )}{s^w_p(1/\omega _0,\vartheta _0)}\nm
{h_1}{M^2_{(\vartheta _1)}}\nm {h_2}{M^2_{(\vartheta _2)}}.
$$
Furthermore, if $h_1=h_2$ and $b^w(x,D)\ge 0$, where $b=a(\sqrt 2\cdo
)$, then $\tp _{h_1,h_2}(a)\ge 0$.
\end{thm}

\par

\begin{proof}
Since $W_{h_2,h_1}\in s_1^w(1/\vartheta _1,\vartheta _2)$, the result
is an immediate consequence of \eqref{ToepWeyl} and Theorem
\ref{thm3.3weyl}.
\end{proof}

\end{document}